\documentclass[11pt]{article}
\usepackage{amsmath} \usepackage{amssymb} \usepackage{latexsym} \usepackage{enumitem} \usepackage{mathrsfs} \usepackage{comment}
\usepackage{color}
\usepackage{accents}
\usepackage{bbm}

\usepackage[colorlinks=true,urlcolor=blue,
citecolor=red,linkcolor=blue,linktocpage,pdfpagelabels,bookmarksnumbered,bookmarksopen]{hyperref}

\newtheorem{assumption}{Assumption}

\def\qed{ \ \vrule width.2cm height.2cm depth0cm\smallskip}
\newenvironment{proof}{\noindent {\bf Proof.\/}}{$\qed$\vskip 0.1in}

\newcommand{\la}{\langle}
\newcommand{\ra}{\rangle}

\newcommand{\ol}{\overline}
\newcommand{\ul}{\underline}

\newcommand{\ba}{\begin{array}}
\newcommand{\ea}{\end{array}}
\newcommand{\be}{\begin{equation}}
\newcommand{\ee}{\end{equation}}
\newcommand{\bea}{\begin{eqnarray}}
\newcommand{\eea}{\end{eqnarray}}
\newcommand{\beaa}{\begin{eqnarray*}}
\newcommand{\eeaa}{\end{eqnarray*}}

\def\dbE{\mathbb{E}}
\def\dbF{\mathbb{F}}

\def\dbL{\mathbb{L}}

\def\dbP{\mathbb{P}}
\def\dbR{\mathbb{R}}

\def\a{\alpha}
\def\b{\beta}
\def\g{\gamma}

\def\e{\varepsilon}

\def\l{\lambda}

\def\f{\varphi}
\def\th{\theta}
\def\o{\omega}
\def\h{\widehat}
\def\G{\Gamma}

\def\O{\Omega}
\def\cA{{\cal A}}

\def\cC{{\cal C}}

\def\cF{{\cal F}}
\def\cG{{\cal G}}

\def\cL{{\cal L}}
\def\cM{{\cal M}}

\def\cP{{\cal P}}

\def\cX{{\cal X}}
\def\cY{{\cal Y}}
\def\cZ{{\cal Z}}

\def\no{\noindent}

\def\ms{\medskip}

\def\q{\quad}
\def\qq{\qquad}

\def\pa{\partial}
\def\cd{\cdot}
\def\cds{\cdots}

\def\td{\nabla}

\def\tr{\hbox{\rm tr}}

\def\qed{ \hfill \vrule width.25cm height.25cm depth0cm\smallskip}

\newcommand{\basa}{\begin{assumption}}
\newcommand{\easa}{\end{assumption}}

\newcommand{\bas}{\begin{assum}}
\newcommand{\eas}{\end{assum}}

\def\limsup{\mathop{\overline{\rm lim}}}
\def\liminf{\mathop{\underline{\rm lim}}}

\def\pa{\partial}
\def\h{\widehat}

 \def\cd{\cdot}
\def\cds{\cdots}

\def\tr{\hbox{\rm tr$\,$}}

\def\dis{\displaystyle}

\def\wh{\widehat}

\def\1{{\bf 1}}

\def\:{\!:\!}
\def\reff{\eqref}
\def \proof{{\noindent \bf Proof.\quad}}

at 9pt

\definecolor{alp}{rgb}{0.0, 0.5, 0.0}

\newtheorem{thm}{Theorem}[section]
\newtheorem{lem}[thm]{Lemma}
\newtheorem{cor}[thm]{Corollary}
\newtheorem{prop}[thm]{Proposition}
\newtheorem{rem}[thm]{Remark}

\newtheorem{defn}[thm]{Definition}
\newtheorem{assum}[thm]{Assumption}

\begin{document}

\title{\bf  Minimal solutions of master equations for extended mean field games} 
\author{Chenchen Mou\thanks{\noindent  Department of Mathematics,
City University of Hong Kong. E-mail: \href{mailto:chencmou@cityu.edu.hk}{chencmou@cityu.edu.hk}. This author is supported in part by Hong Kong RGC Grants ECS 21302521, GRF 11311422 and GRF 11303223.} ~ and ~ Jianfeng Zhang\thanks{\noindent  Department of Mathematics, 
University of Southern California. E-mail:
\href{mailto:jianfenz@usc.edu}{jianfenz@usc.edu}. This author is supported in part by NSF grants DMS-1908665 and DMS-2205972.
}  
}
\date{}
\maketitle

\begin{abstract} 
In an extended mean field game the vector field governing the flow of the population can be different from that of the individual player at some mean field equilibrium. This new class strictly includes the standard mean field games. It is well known that, without any monotonicity conditions, mean field games typically contain multiple mean field equilibria and the wellposedness of their corresponding master equations fails. In this paper, a partial order for the set of probability measure flows is proposed to compare different mean field equilibria. The minimal and maximal mean field equilibria under this partial order are constructed  and satisfy the flow property. The corresponding value functions, however, are in general discontinuous.  We thus introduce a notion of weak-viscosity solutions for the master equation and verify that the value functions are indeed weak-viscosity solutions. Moreover, a comparison principle for weak-viscosity semi-solutions is established and thus these two value functions serve as the minimal and maximal weak-viscosity solutions in appropriate sense.  
 In particular, when these two value functions coincide, the value function becomes the unique weak-viscosity solution to the master equation. The novelties of the work persist even when restricted to the standard mean field games.
\end{abstract}

\no{\bf Keywords.}  Mean field games, master equation, weak solution, viscosity solution, comparison principle.

\ms
\no{\it 2020 AMS Mathematics subject classification:}  35Q89, 49N80, 35D40, 60H30, 91A16, 93E20

\vfill\eject


\section{Introduction}
\label{sect-Introduction}
\setcounter{equation}{0} 
In this paper we consider the  following extended mean field game system:
 given $\mu\in \cP_2(\dbR^d)$,
\bea
\label{MFGsystem}
\left.\ba{c}
\dis \pa_t \nu(t,x) - {1\over 2} \tr\big(\pa_{xx} \nu(t,x)) + div(\nu(t,x) \wh b(x, \pa_x v(t,x), \nu_t)\big) =0,\q \nu_0 = \mu;\ms\\
\dis \pa_t v(t,x) + {1\over 2} \tr(\pa_{xx} v(t,x)) + H(x, \pa_x v(t,x), \nu_t) =0,\q v(T,x) = G(x, \nu_T).
\ea\right.
\eea
The master equation, see \reff{master} below, is to characterize its decoupling field $V$ in the sense that
\bea
\label{vV}
v(t,x) = V(t, x, \nu_t).
\eea
The standard mean field game and its master equation correspond to the special case:
\bea
\label{whb=b}
\wh b(x, p, \mu) = \pa_p H(x, p, \mu).
\eea

Initiated independently by Caines-Huang-Malham\'e \cite{HCM06} and Lasry-Lions \cite{LL07a}, mean field games (MFGs, for short) have received very strong attention and is by now a well-established theory for the study of the asymptotic behavior of stochastic differential games with a large number of players interacting in certain symmetric way. We refer to the monographs Carmona-Delarue \cite{CD1,CD2} and the lecture note Cardaliaguet-Porretta \cite{CP} for a complete introduction of recent progresses on the subject. 

Extended MFGs were first introduced by Lions-Souganidis \cite{LionsSouganidis} to study a more general class of MFGs where the vector field governing the flow of the population might be different from that of the individual player at some mean field equilibrium (MFE, for short). 
Their motivation comes from two folds. Firstly,  the  homogenization limit of a class of oscillatory classical MFGs is in general not a classical MFG but an extended MFG. Secondly, extended MFGs arise naturally in the optimal transportation-type control problems. More precisely, the Euler-Lagrange systems of optimal transportation-type control problems are in general not of the classical MFG type but of the extended MFG type. A new and meaningful monotonicity condition was proposed in \cite{LionsSouganidis} to study the wellposedness of extended MFG systems, and their wellpoedness results were further extended in Mun\~oz \cite{Munoz}. In particular,  the proposed monotonicity condition ensures the uniqueness of MFE of extended MFGs.

It should be noted that \cite{LionsSouganidis, Munoz} consider extended MFG systems with local coupling, that is, the data $G, H, \wh b$ depend on $\nu(t,x)$, rather than $\nu_t$. We instead study extended MFGs with nonlocal coupling, as in \reff{MFGsystem}, via the  master equation \eqref{master}. Our motivation for studying such extended MFGs comes from the study of MFGs with a major player. These games consist of a major player and infinite many homogeneous minor players where the major player can have a significant impact on the minor players while all the minor players as a whole can have an impact on the major player. In this case, the value function of the major player will take the form
\bea
\label{Vextended}
V_0(t, X^0_t, \cL_{X_t|\cF^{X^0}_t}),
\eea
where $X^0$ and $X$ stand for the major player's state and the representative minor player's state, respectively. In particular, the measure variable $\cL_{X_t|\cF^{X^0}_t}$ is not the law of the major player's state $X^0_t$. This is exactly in the spirit of the extended MFG. The local (in time) wellposedness of the MFG systems for MFGs with a major player has been established in Cardaliaguet-Cirant-Porretta \cite{CCP}.  Its global wellposedness has not been studied in the literature, to the best of our knowledge, and we shall address it in our accompanying paper \cite{MZ-majorMFG}.

In the literature of standard MFGs, the global wellposedness of master equations requires the uniqueness of MFE, typically under certain monotonicity conditions. 
See, e.g., Bertucci \cite{B1}, Bertucci-Cecchin \cite{BC}, Cardaliaguet-Delarue-Lasry-Lions \cite{CDLL}, Cardaliaguet-Souganidis \cite{CarSou}, Carmona-Delarue \cite{CD2},  Chassagneux-Crisan-Delarue \cite{CCD},  Lions \cite{Lions}, Mou-Zhang \cite{MZ2},  for the well-known Lasry-Lions monotonicity condition; Ahuja \cite{Ahuja}, Bensoussan-Graber-Yam \cite{BGY1,BGY2}, Gangbo-Meszaros \cite{GM}, Gangbo-Meszaros-Mou-Zhang \cite{GMMZ} for the displacement monotonicity condition; and Mou-Zhang \cite{MZ3} for the anti-monotonicity condition. We emphasize that, all these monotonicity conditions require the measure variable to be the law of the state process, and thus fail automatically for value functions in the form \reff{Vextended}. The works Graber-Meszaros \cite{GraberM1, GraberM2} proposed a new type of monotonicity condition, which does not have this constraint. We should mention the very recent work Bertucci-Lasry-Lions \cite{BLL} concerning master equations for extended MFGs with nonlocal coupling as in the present paper. It shows that the master equation admits at most one global solution which is Lipschitz continuous in the measure variable. However, the existence of such a solution requires additional structural conditions and remains open. Moreover, there are studies on master equations for finite state extended MFGs, see e.g. Bertucci \cite{B2} and Bertucci-Lasry-Lions \cite{BLL1,BLL2}. 
We shall investigate the existence of  global classical solutions  of master equations for extended MFGs in the accompanying paper \cite{MZ-majorMFG}.  

In this paper we focus on extended MFGs and their master equations, with possibly multiple MFEs. Our main idea is to introduce a partial order $\preceq$ for the set of probability measure flows, in the spirit of stochastic dominance. This allows us to compare different MFEs, and we shall construct the minimal/maximal MFE for extended MFGs under this partial order, following the Knaster–Tarski fixed point theorem. To be precise, we shall construct MFEs $\ul \nu$ and $\ol \nu$ such that: 
\bea
\label{ulolnu}
\ul \nu \preceq \nu^* \preceq \ol \nu,\q\mbox{for all MFE $\nu^*$}.
\eea
For this purpose, we shall assume the data $G, H, \wh b$ are monotone in $\mu$ under the partial order $\preceq$. We emphasize that this type of monotonicity under $\preceq$ has a completely different nature from the various monotonicity conditions mentioned in the previous paragraph. Our approach is strongly inspired by Dianetti-Ferrari-Fischer-Nendel \cite{DFFN1,DFFN2} and Dianetti \cite{Dianetti} which obtained \reff{ulolnu} under the same partial order for standard MFGs. A similar idea has also been applied previously to investigate MFGs of optimal stopping, see Carmona-Delarue-Lacker \cite{CDL} and Bertucci \cite{B0}.

We next establish the flow property of the minimal/maximal MFEs, which is crucial for studying the dynamic value function and the master equation. That is, let $\ul \nu^{t, \mu}$ denote the minimal MFE for the extended MFG on $[t, T]$ with initial distribution $\mu$. Then, for any $t_0<t_1$,
\bea
\label{flow0}
\ul \nu^{t_0,\mu}_t = \ul \nu^{t_1, \ul \nu^{t_0,\mu}_{t_1}}_t,\q t\ge t_1.
\eea
This implies the following value function is time consistent:
\bea
\label{ulV0}
\ul V(t_0, x, \mu) = v(t_0, x),\q\mbox{where $v$ solves the backward PDE in \reff{MFGsystem} with $\nu = \ul\nu^{t_0,\mu}$}.
\eea
This function $\ul V$ is smooth in $x$, but is typically discontinuous in $(t, \mu)$, as we will see in Section \ref{sect-eg} below. So a classical solution theory for the master equation is not viable under our conditions. 

We thus turn to weak solutions, by adapting the notion of weak-viscosity solution proposed in our previous paper \cite{MZ2}. We shall show that, by introducing $\ol V$ associated to the maximal MFE, both $\ul V$ and $\ol V$ are weak-viscosity solutions of the master equation \reff{master}. Moreover, for any weak-viscosity solution $V$, the spatial derivative $\pa_x V$ always stays between $\pa_x \ul V$ and $\pa_x \ol V$ component wise. In this sense, $\ul V$ and $\ol V$ can be viewed as the minimal and maximal weak-viscosity solution of the master equation.  In particular, the weak-viscosity solution is unique if and only if $\ul V = \ol V$. We would like to note that, the very recent work Lions-Seeger \cite{LionsSeeger} has used the same approach  to establish the global well-posedness for linear and nonlinear finite dimensional transport equations with coordinate-wise increasing velocity fields, and the theory has also been applied to study MFGs in a finite state space.

We note that our consideration of $\ul \nu$ and $\ol \nu$ can be viewed as a special selection of MFEs. In the literature there have been  other selection criteria for standard MFGs with multiple MFEs, see e.g. Delarure-Foguen Tchuendom \cite{DF}, Cecchin-Dai Pra-Fisher-Pelino \cite{CDFP}, Cecchin-Delaure \cite{CecchinDelarue,CecchinDelarue1}. In \cite{DF}, three methods of selection, including the minimal cost, zero noise limit, $N$-player limit selections, are considered for the linear quadratic MFGs. In particular, in this case the master equation is reduced to a one dimensional PDE and the MFE selected by the last two methods provides an entropy solution to this PDE.
Similar results have been obtained for two-state MFGs in \cite{CDFP}. In  \cite{CecchinDelarue,CecchinDelarue1} the authors established the global wellposedness of master equations for potential MFGs with multiple MFEs. The potential game structure allows to link the MFG to a mean field control problem in the sense that the selected MFE for the MFG is an optimal strategy for the control  problem. We would also like to mention that Iseri-Zhang \cite{IZ} takes a different approach by investigating the set value of MFGs, namely the set of game values over all MFEs, which satisfies the dynamic programming principle. Again our $\ul V$ and $\ol V$ can be viewed as the minimal and maximal (in terms of $\pa_x V$ instead of $V$)  elements  of the set value. 

The rest of the paper is organized as follows. In Section \ref{sec:setting} we introduce the problem, the main results, and the assumptions. In Section \ref{sec:HJB} we investigate the backward PDE in \reff{MFGsystem} for given $\nu$. In Section \ref{sec:minMFE} we construct the minimal MFE for the extended MFG. In Section \ref{sec:value}
 we study the basic properties of the value function $\ul V$.   In Section \ref{sect-minmax} we establish the weak-viscosity solution theory. In Section \ref{sect-extension} we present the results concerning the maximal MFE and its corresponding value function $\ol V$;  the results under an alternative set of monotonicity condition under the partial order; as well as the extension of the current results to extended MFGs with a common noise. Finally in Section \ref{sect-eg} we solve an example explicitly, which in particular shows that $\ul V$ is discontinuous in $(t, \mu)$.

\section{The setting and the main results} 
\label{sec:setting}
\setcounter{equation}{0}

Throughout the paper, we fix a  finite time horizon $[0, T]$ and a filtered probability space $(\O,\cF, \dbF, \dbP)$, on which is defined a $d$-dimensional Brownian motion $B$. For any $p\ge 1$, let $\cP_p(\dbR^d)$ denote the set of probability measures on $\dbR^d$ with finite $p$-th moment, equipped with the $p$-Wasserstein distance $W_p$. We assume $\cF_0$ is rich enough to support any $\mu\in \cP_2(\dbR^d)$, and $\cF_t := \cF_0\vee\cF_t^B$.  For any $p\ge 1$, $\cG\subset \cF$, and $\mu\in \cP_p(\dbR^d)$, denote by $\dbL^p(\cG)$ the set of $\cG$-measurable and $p$-integrable random variables $\xi$; and $\dbL^p(\cG;\mu)$ the set of those $\xi\in \dbL^p(\cG)$ with $\cL_\xi=\mu$. For any $t_0\in [0, T]$, denote $B_t^{t_0}:=B_t-B_{t_0}$, $t\in [t_0,T]$, and $\dbF^{t_0}:=\{\cF_t\}_{t_0\le t\le T}$. Moreover, we denote ${\bf 0} := (0,\cds, 0)$ and ${\bf 1} := (1,\cds, 1)$ with appropriate dimensions.

\subsection{The extended mean field game}
\label{sec:EMFG}
\setcounter{equation}{0}
First, given $t_0\in [0, T]$ and $\nu\in C([t_0,T]; \cP_2(\dbR^d))$, consider the following parabolic PDE on $[t_0, T]$:
\bea
\label{HJB}
 \pa_t v(\nu; t,x)+\frac{1}{2}\tr(\pa_{xx}v(\nu; t,x))+H(x,\pa_xv(\nu; t,x),\nu_t)=0,
\qq v(\nu; T,x)=G(x,\nu_T).
\eea
Under certain technical conditions on $H, G$ as we will specify later, the above PDE has a unique classical solution $v(\nu; \cd,\cd)$. Next, given $\xi\in \dbL^2(\cF_{t_0})$, consider the following SDE on $[t_0, T]$:
\bea
\label{Xb}
X^{t_0, \xi, \nu}_t = \xi + \int_{t_0}^t \wh b\big(X_s^{t_0,\xi,\nu}, \pa_xv(\nu; s,X_s^{t_0,\xi,\nu}),\nu_s\big)ds+B_t^{t_0}.
\eea
It is clear that the mapping $\xi \mapsto \cL_{X^{t_0, \xi, \nu}}$ is law invariant.  We then define the Nash field $\Phi$ for the extended MFG as follows:  for any $(t_0, \mu)\in [0, T]\times \cP_2(\dbR^d)$ and $\xi\in \dbL^2(\cF_{t_0};\mu)$, 
\bea
\label{NashField}
\dis \Phi(t_0, \mu, \nu) := \{\cL_{X^{t_0, \xi, \nu}_t}\}_{t_0\le t\le T},\q\forall \nu\in C([t_0, T]; \cP_2(\dbR^d)).
\eea

\begin{defn}
\label{defn-MFE}
For any $(t_0,\mu)\in [0,T]\times\mathcal{P}_2(\mathbb R^d)$,   we say $\nu^*\in C([t_0,T]; \cP_2(\dbR^d))$ is a mean field equilibrium (MFE) at $(t_0,\mu)$  if it is a fixed point of the Nash field $\Phi(t_0, \mu, \cd)$:
\bea\label{MFE}
\Phi(t_0, \mu,\nu^*) = \nu^*.
\eea
\end{defn}

\begin{rem}
\label{rem-MFE}
(i) The typical case is that $H$ is a Hamiltonian and thus \reff{HJB} is the HJB equation:
\bea
\label{eq:H}
H(x,p,\mu):=\inf_{a\in\dbR}h(x,p,\mu,a),\quad\mbox{where}\q  h(x,p,\mu,a):=p\cd b_0(x,a,\mu)+f(x,a,\mu).
\eea
In this case, as in the standard theory we have a representation formula for $v$: 
\bea
\label{vnu}
\left.\ba{lll}
\dis X^{0, \nu; t_0, x, \a}_t = x+\int_{t_0}^{t}b_0(X^{0,\nu; t_0, x, \a}_s,\alpha(s, X^{0,\nu; t_0, x, \a}_s), \nu_s)ds+B_t^{t_0};\\
\dis J( \nu; t_0, x, \a) := \dbE\Big[g(X^{0, \nu; t_0, x, \a}_T, \nu_T) + \int_{t_0}^T f(X^{0,\nu; t_0, x, \a}_s,\alpha(s, X^{0,\nu; t_0, x, \a}_s), \nu_s)ds\Big];\\
\dis v(\nu; t_0, x) := \inf_{\a\in \cA_{t_0}} J( \nu; t_0, x, \a) 
\ea\right.
\eea
where $\cA_{t_0}$ denotes the appropriate set of admissible controls $\a:[t_0,T]\times\dbR^d\to\dbR^d$. \\
(ii) In the case in which the Hamiltonian $H$ has a  minimizer $a^* = \phi(x, p, \mu)$, namely
\bea
\label{I}
H(x,p,\mu) = h(x,p, \mu, \phi(x, p, \mu)).
\eea
By \reff{eq:H}  one can  easily check that
\bea
\label{Hp}
b_0(x,\phi(x,p,\mu),\mu)=\pa_pH(x,p,\mu), \q  f(x,\phi(x,p,\mu),\mu)= H(x, p, \mu) - p\cd \pa_pH(x,p,\mu).
\eea
(iii)  Assuming \reff{I} holds true, one typical case of $\wh b$ is: for some appropriate function $b$,
\beaa
\wh b(x,p,\mu) = b(x, \phi(x, p, \mu), \mu).
\eeaa
When $b=b_0$ or $\wh b(x, p, \mu) = \pa_p H(x, p, \mu)$, the extended MFG becomes a standard MFG. 
\end{rem}

\subsection{The master equation}

When there is a unique MFE for each $(t_0,\mu)\in [0,T]\times\cP_2(\dbR^d)$, denoted as $(\a^*(t_0,\mu;\cd),\nu^*(t_0,\mu))$. Then the game problem leads to the following value function:
\begin{equation}\label{eq:candidate}
V(t_0,x,\mu):=J(\nu^*(t_0,\mu); t_0,x,\a^*(t_0,\mu;\cd))\quad\text{for any $x\in\dbR^d$.}
\end{equation}
Recall the extended MFG \reff{HJB}, \reff{Xb}, \reff{NashField}, and \reff{MFE}.  In light of  \reff{vnu} and \reff{Hp} we introduce the following FBSDE system (the system does not require the structure in Remark \ref{rem-MFE} (i) though): 
\bea
\label{Y*}
\left.\ba{c}
\dis X_t^{0,*}=x+\int_{t_0}^{t}\pa_pH(X_s^{0,*},\pa_xV(s,X_s^{0,*},\nu_s^*),\nu_s^*)ds+B_t^{t_0};\\
\dis X_t^*=\xi+\int_{t_0}^t\wh b(X_s^{*},\pa_xV(s,X_s^{*},\nu_s^*), \nu_s^*)ds+B_t^{t_0};\\
\dis Y_t^{*}= G(X^{0,*}_T, \nu_T^*) - \int_t^T Z^*_s dB_s \\
\dis  +\int_{t}^T \Big[ H(\cd) - \pa_x V(s, X^{0,*}_s, \nu^*_s) \cd \pa_p H(\cd)\Big]\Big(X_s^{0,*}, \pa_x V(s, X^{0,*}_s,\nu^*_s), \nu^*_s\Big)ds;\\
\dis \mbox{where}\q \nu^*_t:= \cL_{X^*_t}.
\ea\right.
\eea
In particular, we have
\bea
\label{YVX}
Y^*_t = V\big(t, X^{0,*}_t, \nu_t^*\big) = V\big(t, X^{0,*}_t, \cL_{X_t^*}\big).
\eea
By applying the It\^{o}'s formula (c.f. \cite{BLPR,CCD}) and comparing it with \eqref{Y*}, we derive the master equation: 
\bea
\label{master}
\left.\ba{c}
\dis \pa_t V + \frac{1}{2} \tr(\pa_{xx} V) +  H(x,\partial_x V,\mu) + \cM V =0, \q V(T,x,\mu) = G(x,\mu),\q \mbox{where}\\
\dis \cM V(t,x,\mu)  := \tr\Big( \int_{\dbR^d}\Big[\frac{1}{2} \pa_{\tilde x} \pa_\mu V(t,x, \mu, \tilde x)  + \pa_\mu V(t, x, \mu, \tilde x)\wh b^\top(\tilde x,\pa_x V(t, \tilde x, \mu),\mu)\Big] \mu(d\tilde x)\Big).
\ea\right.
\eea
Note that we may alternatively view $V$ as the decoupling field of the following FBSDE system: 
\bea
\label{cY*}
\left.\ba{c}
\dis \cX_t^{0,*}=x+B_t^{t_0};\\
\dis \cX_t^*=\xi+\int_{t_0}^t\wh b(X_s^{*},\pa_xV(s,X_s^{*},\nu_s^*), \nu_s^*)ds+B_t^{t_0},\q \mbox{where}\q \nu^*_t:= \cL_{\cX^*_t};\\
\dis \cY_t^{*}= G(\cX^{0,*}_T, \nu_T^*)   +\int_{t}^T H\big(\cX_s^{0,*}, \pa_x V(s, \cX^{0,*}_s,\nu^*_s), \nu^*_s\big)ds- \int_t^T \cZ^*_s dB_s;\\
\dis \mbox{in the sense}\q \cY^*_t = V(t, \cX^{0,*}_t, \nu^*_t).
\ea\right.
\eea
Moreover, $V$ also serves as the decoupling field of the extended MFG system, see \reff{MFGsystem} and \reff{vV}.

The main feature here is that the measure variable $\nu^*_t$ in \reff{YVX} is the law of $X^*_t$, rather than that of $X^{0,*}_t$. Consequently, the $\cM V$ above involves the term $\pa_\mu V \wh b^\top$, instead of $\pa_\mu V b_0^\top= \pa_\mu V \pa_p H^\top$ as in the standard master equations. This feature appears naturally in MFG with a major player, which is the main motivation of this paper and will be the subject of our accompanying paper \cite{MZ-majorMFG}.
We also refer to  \cite{LionsSouganidis} for more applications of extended MFGs.

However, in general there could be multiple MFEs, which lead to multivalued functions. Our goal in this paper is to construct the minimal/maximal MFE and to verify that their value functions satisfy the master equation, in the sense of weak-viscosity solutions introduced in \cite{MZ2}.

\subsection{The main results}
The main results of this paper build on the following partial order $\preceq$ (or alternatively $\succeq$).
\begin{defn}
\label{defn-order} For a generic dimension $n$ and for $i=1, 2$,\\
(i) for any $x^i = (x^i_1,\cds, x^i_n)\in\dbR^n$, we say that $x^1\preceq x^2$ if $x^1_j\leq x^2_j$ for all $j=1,\cds, n$;\\
(ii) for any $\mu_i\in\cP_2(\mathbb R^n)$, we say that $\mu_1\preceq \mu_2$ if there exist $\xi^i\in\dbL^2(\mathcal{F}_0;\mu_i)$  s.t. $\xi^1\preceq \xi^2$ $\dbP$-a.s.;\\
(iii) for any $\nu^i\in C([t_0,T];\cP_2(\dbR^n))$, we say that $\nu^1\preceq\nu^2$ if $\nu^1_t\preceq\nu^2_t$ for all $t\in [t_0,T]$.
\end{defn}
We note that $\mu_1\preceq \mu_2$ is equivalent to the stochastic dominance. We say $x^1 \succeq x^2$ if $x^2 \preceq x^1$, and  a function $\f: \dbR^n \to \dbR^m$ is increasing (resp. decreasing) if $\f(x^1) \preceq \f(x^2)$ whenever $x^1 \preceq \mbox{(resp. $\succeq$)}~ x^2$. Similarly we define the monotonicity of functions on  $\cP_2(\dbR^d)$ and $C([t_0,T];\cP_2(\dbR^d))$.

We first have the following simple proposition. 
\begin{prop}\label{prop:inc}
Assume $\f\in \cC^1(\cP_2(\mathbb R^d))$, namely it has a continuous Lions derivative $\pa_\mu \f$. Then $\f$  is increasing if and only if  $\pa_{\mu}\f(\mu, x) \succeq {\bf 0}$ for all $(\mu, x)\in \cP_2(\dbR^d)\times\dbR^d$.
\end{prop}

\proof We first prove the if part. Assume $\pa_{\mu}\f \succeq {\bf 0}$. Let $\mu_1,\mu_2\in\cP_2(\dbR^d)$ be such that $\mu_1\preceq\mu_2$, i.e. there exist $\xi^i\in\dbL^2(\cF_0;\mu_i)$, $i=1,2$, such that $\xi^1\preceq \xi^2$ $\dbP$-a.s. Then
\[
\f(\mu_2)-\f(\mu_1)=\int_{0}^1\dbE\Big[\pa_\mu \f\big(\cL_{\xi^1+\theta(\xi^2-\xi^1)},\xi^1+\theta(\xi^2-\xi^1)\big)\cd(\xi^2-\xi^1)\Big]d\th\geq 0.
\]

We next prove the only if part. Assume $\f$ is increasing. For any $\mu\in\cP_2(\dbR^d)$,  $\xi\in \dbL^2(\cF_0;\mu)$, and $\eta\in \dbL^2(\cF_0)$ such that $\eta \succeq {\bf 0}$,  we have
\[
0\leq \lim_{\e\downarrow 0}\frac{\f(\cL_{\xi+\e\eta})-\f(\mu)}{\e}=\dbE\Big[\pa_\mu \f(\mu,\xi)\cd\eta\Big].
\]
By the arbitrariness of  $\eta \succeq {\bf 0}$, this implies that $\pa_{\mu}\f(\mu,\xi)\succeq {\bf 0}$, $\dbP$-a.s. That is $\pa_{\mu}\f(\mu,x)\succeq {\bf 0}$, for $\mu$-a.e. $x$. Since $\pa_\mu \f$ is continuous, we see that $\pa_{\mu}\f(\mu, x) \succeq {\bf 0}$ for all $(\mu, x)$.
\qed

\begin{rem}
\label{rem-mon}
As we saw in \cite{GMMZ}, a smooth function $U$ on $\dbR^d\times \cP_2(\dbR^d)$ satisfies the Lasry-Lions monotonicity condition if and only if: for any $\mu\in \cP_2(\dbR^d)$, $\xi\in \dbL^2(\cF_0; \mu)$, $\eta\in \dbL^2(\cF_0)$, 
\bea
\label{LLmon}
\dbE\Big[\la\pa_{x\mu} U(\xi, \mu, \tilde \xi) \eta,~ \tilde \eta\ra\Big]\ge 0.
\eea
We note that \reff{LLmon} is always under expectation, while in Proposition \ref{prop:inc} we require $ \pa_{\mu}\f(\mu, x) \succeq {\bf 0} $ pointwisely. In this sense we are considering pointwise monotonicity in this paper. We shall remark that \reff{LLmon} and the pointwise monotonicity of $\pa_x U(x, \cd)$ do not imply each other. 
\end{rem}

Our main results consist of two parts, under the conditions specified in the next subsection. 

\begin{itemize}
\item{} First, given $(t, \mu)\in [0,T]\times\cP_2(\dbR^d)$, we will construct the minimal MFE $\ul \nu^{t, \mu}$ and the maximal MFE  $\ol \nu^{t, \mu}$ at $(t, \mu)$, in the sense that for any other MFE $\nu^*$ at $(t, \mu)$ it holds:
\beaa
\ul \nu^{t, \mu} \preceq \nu^* \preceq  \ol \nu^{t, \mu}.
\eeaa

\item{} Next, we define the dynamic value functions 
\beaa
\ul V(t,x,\mu) := v(\ul \nu^{t, \mu}; t, x),\q \ol V(t,x,\mu)  := v(\ol \nu^{t, \mu}; t, x).
\eeaa
We shall show that they are weak-viscosity solutions of the master equation \reff{master} such that $\pa_x \ul V$ and $\pa_x \ol V$ satisfy certain minimal/maximal property.
\end{itemize}

\no Since the analyses are similar, in the paper we will focus only on $\ul \nu^{t, \mu}$ and $\ul V(t,x,\mu)$, and we will present the results concerning  $\ol \nu^{t, \mu}$ and $\ol V(t,x,\mu)$ in Section \ref{sect-max} below. 
\subsection{The  assumptions}
\label{sec:Assum}

We first introduce some technical assumptions on the coefficients, which are more or less standard in the literature. Denote, for any $R>0$,
\begin{equation}\label{eq:BR}
O_R:=\{p\in\dbR^d\,:\,|p|< R\},\q \forall R>0.
\end{equation}

\begin{assum}\label{assum-GH}
(i)  $G\in C^0(\dbR^d\times\cP_2(\dbR^d))$ and $H\in C^0(\dbR^{2d}\times\cP_2(\dbR^d))$ are functions satisfying $G(\cd,\mu)\in C^2(\dbR^d)$ and  $H(\cd,\cd,\mu)\in C^2(\dbR^d\times\dbR^d)$ for each $\mu\in\cP_2(\dbR^d)$;\\
(ii) there exist constants $L^G_0$, $L^H_0$, and $L^{H}(R)$ for each $R>0$,  such that 
\beaa
&\dis |\pa_xG(x,\mu)|,|\pa_{xx}G(x,\mu)|\leq L_0^G, \q\mbox{and}\q  |\pa_xH(x,p,\mu)|\leq L^H_0[1+|p|],\quad\text{for all}~(x, p, \mu);\\
&\dis |\pa_pH|,|\pa_{xx}H|,|\pa_{xp}H|,|\pa_{pp}H|\leq L^{H}(R)\quad \text{on $\dbR^d\times O_R\times\cP_2(\dbR^d)$;}
\eeaa
(iii) for each $R>0$ and any compact set $K\subset \cP_2(\dbR^d)$, $\pa_xG,\pa_{xx}G$ are uniformly continuous in $(x,\mu)$ on $\dbR^d\times K$, and $\pa_{x}H$, $\pa_{p}H$, $\pa_{xx}H$, $\pa_{xp}H$, $\pa_{pp}H$ are uniformly continuous in $(x, p, \mu)$ on $\dbR^d\times O_R\times K$.
\end{assum}

\begin{assum}\label{assum-b}
Assume that $\wh b(\cd,\cd, \mu)\in C^1(\dbR^d\times \dbR^d)$ for each $\mu\in\cP_2(\dbR^d)$, and for each $R>0$ and any compact set $K\subset\cP_2(\dbR^d)$,  $\wh b$,  $\pa_x \wh b, \pa_p \wh b$ are bounded with bound $L^{\wh b}(R)$ and $\wh b$ is uniformly  continuous in $\mu$ on $\dbR^d\times O_R \times K$.
\end{assum}

The following pointwise monotonicity condition under partial order $\preceq$ is crucial.

\begin{assum}\label{assum-mon}
(i) $\pa_{x}G$ is increasing in $(x, \mu)$;\\
(ii) $\pa_x H$ is increasing in $(x,\mu)$, $\pa_pH$ is increasing in $(p, \mu)$, and $\pa_{x_i p_j} H \ge 0$ for all $i\neq j$ (which is slightly weaker than that $\pa_p H$ is increasing in $x$);\\
(iii) $\wh b$ is increasing in $(p, \mu)$ and $\pa_{x_j} \wh b_i \ge 0$ for all $i\neq j$.
\end{assum}

\no Alternatively, we may replace the above assumption with the following monotonicities. 
\begin{assum}\label{assum-mon2}
(i) $\pa_{x}G$ is decreasing in $(x, \mu)$;\\
(ii) $\pa_x H$ is decreasing in $(x,\mu)$, $\pa_pH$ is increasing in $(p,\mu)$, and $\pa_{x_i p_j} H \ge 0$ for all $i\neq j$;\\
(iii) $\wh b$ is decreasing in $p$, increasing in $\mu$, and $\pa_{x_i} \wh b_j \ge 0$ for all $i\neq j$.
\end{assum}

\no In the paper we will focus only on the analyses under Assumption \ref{assum-mon}. The corresponding results under Assumption \ref{assum-mon2} are essentially the same, with obvious changes, so we will present them in Section \ref{sect-decrease} without proofs.

\subsection{Some preliminary comparison results}
In this subsection we present two well known comparison results for multidimensional SDEs and BSDEs, which will play an important role in the paper. The proofs are rather standard, and we refer to \cite{HP} for further discussions on the BSDE case.

\begin{lem}
\label{lem-SDEcomparison}
Consider the following two $n$-dimensional SDE systems: for $k=1,2$,
\bea
\label{barX}
X^{k, i}_t = \xi_k^i + \int_0^t b_k^i (s, X^k_s) ds + B_t^i,\q i=1,\cds, n,
\eea
where $\xi_k^i\in \dbL^2(\cF_0)$ and $b_k^i: [0, T]\times \O\times \dbR^n \to \dbR$ is $\dbF$-progressively measurable.  Assume\\
(i) for $k=1, 2$, $b_k$ is uniformly Lipschitz continuous in $x$ and $\dbE[\int_0^T |b_k(t, 0)|^2dt ] <\infty$; \\
(ii) $b^i_1$ (or $b^i_2$) is increasing in $x_j$ for any $i\neq j$, and  $\xi_1 \preceq \xi_2$ and $b_1\preceq b_2$.\\
Then $X^1_t \preceq X^2_t$, $0\le t\le T$, $\dbP$-a.s.
\end{lem}

\begin{lem}
\label{lem-BSDEcomparison}
Consider the following two $n$-dimensional BSDE systems: for $k=1,2$,
\bea
\label{barYZ}
Y^{k, i}_t = \xi_k^i + \int_t^T f_k^i (s, Y^k_s, Z^{k,i}_s) ds -\int_t^T Z^{k,i}_s \cd d B_s,\q i=1,\cds, n,
\eea
where $\xi_k^i\in \dbL^2(\cF_T)$ and $f_k^i: [0, T]\times \O\times \dbR^n \times \dbR^d\to \dbR^d$ is $\dbF$-progressively measurable.  Assume\\
(i) for $k=1, 2$, $f_k$ is uniformly Lipschitz continuous in $(y,z)$ and $\dbE[\int_0^T |f_k(t, 0,0)|^2dt ] <\infty$; \\
(ii) $f^i_1$ (or $f^i_2$) is increasing in $y_j$ for any $i\neq j$, and  $\xi_1 \preceq \xi_2$ and $f_1\preceq f_2$.\\
Then $Y^1_t \preceq Y^2_t$, $0\le t\le T$, $\dbP$-a.s.
\end{lem}

\section{The PDE \reff{HJB}}
\label{sec:HJB}
\setcounter{equation}{0}
In this section we focus on the properties of the solution $v$ for the PDE  \reff{HJB}. The following lemma is more or less standard. For the sake of completeness, we sketch a proof here. In particular, our probabilistic arguments will remain valid for the common noise case which will be discussed in Section \ref{sect-common} below.

\begin{lem}\label{lem-vreg}
Let Assumption \ref{assum-GH}  hold. \\
(i) For any given $\nu\in C([0,T];\cP_2(\dbR^d))$, the equation \eqref{HJB} admits a unique classical solution $v$, and there exist constants $C_1, C_2>0$, depending on $T$, $d$, $L^G_0$, $L^H_0$, and the function $L^H$, but independent of $\nu$, such that
 \begin{equation}\label{eq:paxvpaxxvbdd}
 |\pa_{x}v|\leq C_1\quad\text{and}\quad |\pa_{xx}v|\le C_2;
 \end{equation}
(ii) for any compact set $K\subset \cP_2(\dbR^d)$, there exists a modulus of continuity function $\rho_K$  such that: for any $\nu,\nu^1,\nu^2\in C([0,T]; \cP_2(\dbR^d))$ satisfying $\nu_t, \nu^1_t, \nu^2_t \in K$ for all $t$,
\bea
\label{eq:abconpaxv}
&\dis |\pa_xv(\nu^1; t,x)-\pa_xv(\nu^2; t,x)|\leq \rho_K\big(\sup_{t\le s\le T}W_2(\nu^1_s,\nu^2_s)\Big);\\
\label{paxvtreg}
&\dis |\pa_xv(\nu; t_1,x)-\pa_xv(\nu; t_2,x)|\leq \rho_K(t_2-t_1),\q \forall 0\le t_1 < t_2 \le T.
\eea
\end{lem}

\proof First it follows from \cite[Proposition 6.1]{GMMZ} that the following function $v(\nu; t, x)$ satisfies \reff{eq:paxvpaxxvbdd}: denoting $X^{t,x}_s := x + B^t_s$, $t\le s\le T$,
\bea
\label{vY}
\left.\ba{c}
\dis v(\nu; t, x) := Y^{t,x,\nu}_t,\q\mbox{where}\\
\dis Y_s^{t,x,\nu}=G(X_T^{t,x},\nu_T)+\int_s^{T}H(X_r^{t,x},Z_r^{t,x,\nu},\nu_r)dr-\int_s^{T}Z_r^{t,x,\nu}\cd dB_r,\q t\le s\le T.
\ea\right.
\eea
In particular, we have
\bea
\label{Zbdd}
|Z_s^{t,x,\nu}| = |\pa_x v(\nu; s, X^{t,x}_s)| \le C_1.
\eea
We note that the assumptions in the statement of \cite[Proposition 6.1]{GMMZ} involve the derivatives of $G$ and $H$ with respect to $\mu$ as well, but they are never used in that proof.

We next  prove \reff{eq:abconpaxv}.  Fix $K$ and let $\rho_K^0$ denote the common modulus of continuity function of $\pa_x G, \pa_{xx} G$ on $\dbR^d\times K$ and that of $\pa_{x}H,\pa_{p}H,\pa_{xx}H,\pa_{xp}H,\pa_{pp}H$ on $\dbR^d\times O_{C_1}\times K$ for the $C_1$ in \reff{eq:paxvpaxxvbdd} or \reff{Zbdd}.  By standard arguments we have
\bea
\label{paxvrep}
\pa_x v(\nu; t,x) = \td_x Y^{t,x,\nu}_t,\q \pa_{xx} v(\nu; t,x) = \td^2_{xx} Y^{t,x,\nu}_t,
\eea
where $\td_x Y^{t,x,\nu}\in \dbR^d$ and $\td^2_{xx} Y^{t,x,\nu}\in\dbR^{d\times d}$ satisfy the following linear BSDEs on $[t, T]$:
\bea
\label{tdY}
&&\dis \left.\ba{lll}
\dis \nabla_{x_i} Y_s^{t,x,\nu}=\pa_{x_i}G(X_T^{t,x},\nu_T) -\int_s^T\nabla_{x_i} Z_r^{t,x,\nu}\cd dB_r\\
\dis \qq\qq\q +\int_s^T[\pa_{x_i}H+\pa_pH\nabla_{x_i} Z_r^{t,x,\nu}](X_r^{t,x}, \td_x Y_r^{t,x,\nu},\nu_r)dr ,\\
\ea\right.\\
\label{td2Y}
&&\dis \left.\ba{lll}\dis \nabla_{x_ix_j} Y_s^{t,x,\nu}=\pa_{x_ix_j}G(X_T^{t,x},\nu_T) -\int_s^T\nabla_{x_ix_j} Z_r^{t,x,\nu}\cd dB_s \\
\dis \q +\int_s^T\Big[\pa_{x_ix_j}H + \sum_{k=1}^d[\pa_{x_i p_k}H\nabla_{x_jx_k}Y_r^{t,x,\nu}+\pa_{x_jp_k}H\nabla_{x_ix_k}Y_r^{t,x,\nu}] \\
\dis\q +\sum_{k,l=1}^d[\nabla_{x_jx_k}Y_r^{t,x,\nu}\pa_{p_kp_l}H\nabla_{x_ix_l}Y_r^{t,x,\nu}]+ \pa_p H \nabla_{x_ix_j} Z_r^{t,x,\nu}\Big](X_r^{t,x}, \td_x Y_r^{t,x,\nu},\nu_r)dr. 
\ea\right.
\eea
Here  we used the fact that $Z^{t,x,\nu}_r = \pa_x v(\nu; r, X^{t,x}_r)=\td_x Y^{t,x,\nu}_r$.
Recall  \reff{Zbdd} again, then we may rewrite \reff{tdY} as: 
\beaa
\left.\ba{lll}
\dis \nabla_{x_i} Y_s^{t,x,\nu}=\pa_{x_i}G(X_T^{t,x},\nu_T)-\int_s^T\nabla_{x_i} Z_r^{t,x,\nu}\cd dB_r\\
\dis\qq +\int_s^T\big[\pa_{x_i}H+\pa_pH (-C_1 \vee \nabla_{x_i} Z_r^{t,x,\nu} \wedge C_1)\big](X_r^{t,x}, \td_x Y_r^{t,x,\nu},\nu_r)dr,
\ea\right.
\eeaa
where the truncation is in the component wise sense. Note that the generator of the above BSDE is Lipschitz continuous. 
Then, by the standard BSDE estimates (cf. \cite[Chapter 4]{Zhang}) 
we can easily obtain \reff{eq:abconpaxv}. Similarly, we can show that $\pa_x v$ and $\pa_{xx} v$ are uniformly continuous in $x$, with a possibly different modulus of continuity function $\rho$.  

Moreover, for any $t_1 < t_2$, note that $\td_x Y^{t_1, x, \nu}_{t_2} = \pa_x v(\nu; t_2, X^{t_1,x}_{t_2})$ and thus, by  \reff{tdY},  
\beaa
&&\dis \pa_x v(\nu; t_1, x) = \td_x Y^{t_1, x, \nu}_{t_1} \\
&&\dis = \pa_x v(\nu; t_2, X^{t_1,x}_{t_2}) + \int_{t_1}^{t_2}[\pa_{x}H+\pa_pH\nabla_{x} Z_r^{t,x,\nu}](X_r^{t,x}, Z_r^{t,x,\nu},\nu_r)dr -\int_{t_1}^{t_2} \nabla_{x} Z_r^{t,x,\nu}\cd dB_r.
\eeaa
Then, noting that  $|\nabla_{x} Z_r^{t,x,\nu}| = |\pa_{xx} v(\nu, r, X^{t,x}_r)| \le C_2$, 
one can easily prove \reff{paxvtreg}, for a possibly different $\rho_K$. Similarly $\pa_x v$ and $\pa_{xx} v$ are also uniformly continuous in $t$. Moreover, since $G$ and $H$ are continuous,  by \reff{vY} one can easily show that $v$ is also continuous in $t$. Then by \reff{vY} clearly $v(\nu; \cd,\cd)$ is the unique classical solution of \reff{HJB}.
\qed

\begin{prop}\label{prop:umon}
Under Assumptions \ref{assum-GH} and \ref{assum-mon} (i)-(ii),  $\pa_x v$ is increasing in $(x, \nu)$.
\end{prop}
\proof First we may rewrite \reff{td2Y} as:  omitting $^{t,x,\nu}$ for notational simplicity,
\bea
\label{td2Y-2}
\left.\ba{lll}
\dis \nabla_{x_ix_j} Y_s=\pa_{x_ix_j}G(X_T,\nu_T) -\int_s^T\nabla_{x_ix_j} Z_r\cd dB_r\\
\dis \q + \int_s^T \Big[f_0\big(r, (\td_{x_kx_l}Y_r)_{(k,l)\neq (i,j)}\big) +\G_r\nabla_{x_ix_j} Y_r + \pa_p H(X_r, \td_x Y_r,\nu_r) \nabla_{x_ix_j} Z_r\Big]dr,\q \mbox{where}\ms\\
\dis \G_r := \Big[\pa_{x_i p_i}H+\pa_{x_jp_j}H+\sum_{l\not =j}\pa_{p_ip_l}H\nabla_{x_ix_l}Y_r+\sum_{k\not=i}\pa_{p_kp_j}H\nabla_{x_jx_k}Y_r\Big](X_r, \td_x Y_r,\nu_r),\\
\dis  f_0\big(r, (y_{k,l})_{(k,l)\neq (i,j)}\big):=\Big[\pa_{x_ix_j}H  + \sum_{k\neq i}\pa_{x_i p_k}Hy_{jk}+\sum_{k\neq j}\pa_{x_jp_k}Hy_{ik} \\
\dis\qq\qq +\sum_{k\not= i,l\not=j}\pa_{p_kp_l}H [(-C_2) \vee y_{j,k} \wedge C_2][(-C_2) \vee y_{i,l} \wedge C_2]\Big](X_r, \td_x Y_r,\nu_r).
\ea\right.
\eea
Here the constant $C_2$ is from \reff{eq:paxvpaxxvbdd} and we used  \reff{paxvrep}.  We may view \reff{td2Y-2} as a $d^2$-dimensional  BSDE system, with index $(i, j)$ and solution $\{(\td_{x_ix_j}Y, \td_{x_ix_j}Z)\}_{(i,j)}$, where $\G$ is viewed as a given coefficient. We next introduce  two  $d^2$-dimensional  BSDE systems, again with index $(i, j)$:
\beaa
\left.\ba{lll}
\dis Y^{1,(i,j)}_s= -\int_s^TZ^{1,(i,j)}_r\cd dB_r + \int_s^T \Big[\G_r Y^{1, (i,j)}_r + \pa_p H(X_r, \td_x Y_r,\nu_r) Z^{1,(i,j)}_r\Big]dr;\ms\\
\dis Y^{2,(i,j)}_s= \pa_{x_ix_j}G(X_T,\nu_T)-\int_s^TZ^{2,(i,j)}_r\cd dB_r \\
\dis\qq + \int_s^T \Big[f_0\big(r, \{(Y^{2, (k,l)}_r)^+\}_{(k,l)\neq (i,j)}\big) +\G_r Y^{2, (i,j)}_r + \pa_p H(X_r, \td_x Y_r,\nu_r) Z^{2,(i,j)}_r\Big]dr.
\ea\right.
\eeaa
By Assumption \ref{assum-mon} (i)-(ii), we have for all $(i,j)$ and $r\in [t,T]$
\beaa
\pa_{x_ix_j}G(X_T,\nu_T) \ge 0,\q f_0\big(r, \{(y_{k,l})^+\}_{(k,l)\neq (i,j)}\big)\ge 0.
\eeaa
Note that $f_0$ is increasing in $\{(y_{k,l})^+\}_{(k,l)\neq (i,j)}$, and it is obvious that $Y^{1, (i,j)}_s\equiv 0$.
Then it follows from Lemma \ref{lem-BSDEcomparison} that 
$Y^{2, (i,j)}_s\ge Y^{1, (i,j)}_s= 0$, and thus $f_0\big(r, \{(Y^{2, (k,l)}_r)^+\}_{(k,l)\neq (i,j)}\big) = f_0\big(r, \{Y^{2, (k,l)}_r\}_{(k,l)\neq (i,j)}\big)$. This  implies that  $\big\{Y^{2, (i,j)}, Z^{2, (i,j)}\big\}_{(i,j)}$ satisfies BSDE system \reff{td2Y-2}. Then $\pa_{x_i x_j} v(\nu; t, x) = \nabla_{x_ix_j}Y_t  = Y^{2, (i,j)}_t\ge 0$. That is, $\pa_x v$ is increasing in $x$. 

Similarly, given $\nu^1, \nu^2\in C([0, T]; \cP_2(\dbR^d))$ such that $\nu^1 \preceq \nu^2$, omit $^{t,x}$ and denote, for $\th\in [0,1]$, 
\beaa
\bar \td_{x_i} Y_s:= \nabla_{x_i} Y_s^{\nu^2}-\nabla_{x_i} Y_s^{\nu^1},\q \bar \td_{x_i} Z_s:= \nabla_{x_i} Z_s^{\nu^2}-\nabla_{x_i} Z_s^{\nu^1},\q \td_x Y^\th_s:= (1-\th) \td_x Y^{\nu^2}_s + \th \td_x Y^{\nu^1}_s.
\eeaa
Note that $\td_{x_i} Z= \big(\td_{x_ix_1} Y,\cds, \td_{x_i x_d} Y)^\top$. By \reff{tdY} we have
\bea
\label{tdY-difference}
 \left.\ba{lll}
\dis \bar\nabla_{x_i} Y_s=[\pa_{x_i}G(X_T,\nu^2_T)-\pa_{x_i}G(X_T,\nu^1_T)] -\int_s^T\bar\nabla_{x_i} Z_r\cd dB_r\\
\dis \q + \int_s^T \Big[\bar \g_r+ \bar f_0\big(r, \{\bar \td_{x_j}Y_r\}_{j\neq i}\big) + \bar\G_r \bar\nabla_{x_i} Y_r +  \pa_pH(X_r, \td_x Y_r^{\nu^1},\nu^1_r) \bar \nabla_{x_i} Z_r\Big]dr,\q\mbox{where}\ms\\
\dis \bar \G_r:=  \int_0^1 \Big[\pa_{x_i p_i} H +  \sum_{k=1}^d\pa_{p_i p_k} H  \td_{x_i x_k} Y^{\nu^2}_r\Big](X_r, \td_x Y^\th, \nu^1)d\th,\\
\dis \bar f_0\big(r, \{y_j\}_{j\neq i}\big):= \sum_{j\neq i} \int_0^1 \Big[\pa_{x_i p_j} H  +  \sum_{k=1}^d   \pa_{p_j p_k} H\td_{x_i x_k} Y^{\nu^2}_r \Big](X_r, \td_x Y^\th, \nu^1)d\th ~ y_j\\
\dis\bar\g_r:= [\pa_{x_i}H(X_r, \td_x Y_r^{\nu^2},\nu^2_r)-\pa_{x_i}H(X_r, \td_x Y_r^{\nu^2},\nu^1_r)] \\
\dis\q+\sum_{k=1}^d [\pa_{p_k}H(X_r, \td_x Y_r^{\nu^2},\nu^2_r) - \pa_{p_k}H(X_r, \td_x Y_r^{\nu^2},\nu^1_r)] \td_{x_i x_k} Y^{\nu^2}_r.
\ea\right.
\eea
Note that $\td_{x_i x_k} Y^{\nu^2}_r = \pa_{x_ix_k} v(\nu^2; r, X_r) \ge 0$. Then, by Assumption \ref{assum-mon} (i)-(ii) we see that $\bar f_0$ is increasing in $\{y_j\}_{j\neq i}$ and, for all $i$ and $r\in [t,T]$,
\beaa
[\pa_{x_i}G(X_T,\nu^2_T)-\pa_{x_i}G(X_T,\nu^1_T)]\ge 0,\q \bar \g_r \ge 0.
\eeaa
Now compare \reff{tdY-difference} with the following $d$-dimensional linear BSDE system:
\bea
\label{tdY-difference2}
\bar Y^i_s= -\int_s^T\bar Z^i_r\cd dB_r+ \int_s^T \Big[ \bar f_0\big(r, \{Y^j_r\}_{j\neq i}\big) + \bar\G_r \bar Y^i_r +  \pa_pH(X_r, \td_x Y_r^{\nu^1},\nu^1_r) \bar Z^i_r\Big]dr.
\eea
%
%
It follows from Lemma \ref{lem-BSDEcomparison} again that $\bar\nabla_{x_i} Y_s \ge \bar Y^i_s$ for all $i$. From \reff{tdY-difference2} it is obvious that $\bar Y^i_s\equiv 0$. Then 
\beaa
\pa_{x_i} v(\nu^2; t, x) - \pa_{x_i} v(\nu^1; t, x) = \bar\nabla_{x_i} Y_t = \bar Y^i_s \ge 0.
\eeaa
 That is, $\pa_x v$ is increasing in $\nu$. 
\qed

\section{The minimal MFE}
\label{sec:minMFE}
\setcounter{equation}{0}
In this section we construct the minimal MFE for the extended MFG. 
We first establish  the pointwise monotonicity of the Nash field $\Phi$.

\begin{thm}\label{thm:increase}
Let Assumptions  \ref{assum-GH}, \ref{assum-b}, and \ref{assum-mon} hold. Then for any $t_0\in [0,T]$, $\Phi(t_0, \cd, \cd)$ is increasing in $(\mu, \nu)$.
\end{thm}
\proof  Let $\mu_1, \mu_2 \in \cP_2(\dbR^d)$ and  $\nu^1, \nu^2 \in C([t_0, T];\cP_2(\dbR^d))$ be such that $\mu_1\preceq \mu_2$, $\nu^1\preceq \nu^2$, and $\xi_1\in \dbL^2(\cF_{t_0}; \mu_1)$, $\xi_2\in \dbL^2(\cF_{t_0}; \mu_2)$ be such that $\xi_1\le \xi_2$. For $k=1,2$, we have
\beaa
X^{t_0, \xi_k, \nu^k}_t = \xi_k + \int_{t_0}^t\wh b(X_s^{t_0,\xi_k,\nu^k}, \pa_xv(\nu^k; s,X_s^{t_0,\xi_k,\nu^k}),\nu^k_s)ds+B^{t_0}_t.
\eeaa
Denote $b_k(s, x):= \wh b(x, \pa_x v(\nu^k; s, x), \nu^k_s)$,  $k=1,2$. By Lemma \ref{lem-vreg} $b_k$ satisfies Lemma \ref{lem-SDEcomparison} (i). Moreover, by Assumption \ref{assum-mon} (iii) and Proposition \ref{prop:umon} we see that $b_1 \preceq b_2$ and 
\beaa
 \pa_{x_j} b_k^i(s, x) = \Big[\pa_{x_j} \wh b^i + \pa_p \wh b^i \cd \pa_{x_j x} v\Big](x, \pa_x v(\nu^k; s, x), \nu^k_s)\ge 0,\q i\neq j.
\eeaa
Since $\xi_1 \preceq \xi_2$, then by Lemma \ref{lem-SDEcomparison} we have $X^{t_0, \xi_1, \nu^1}_t\preceq X^{t_0, \xi_2, \nu^2}_t$, $t_0\le t\le T$, $\dbP$-a.s. This implies that $\Phi(t_0, \mu_1, \nu^1) \preceq \Phi(t_0, \mu_2, \nu^2)$.
\qed

We now construct the minimal MFE by Picard iteration, following the standard procedure in Knaster-Tarski fixed point theorem.  Fix $(t_0, \mu)\in [0, T]\times \cP_2(\dbR^d)$ and $\xi\in \dbL^2(\cF_{t_0}; \mu)$. Recall Assumption \ref{assum-b} and \reff{eq:paxvpaxxvbdd}, we set
\bea
\label{ulX0}
\ul X^{t_0, \xi, 0}_t := \xi - L^{\wh b}(C_1) {\bf 1} + B^{t_0}_t,\q  \ol X^{t_0, \xi, 0}_t := \xi + L^{\wh b}(C_1) {\bf 1} + B^{t_0}_t,
\eea
and, for $n=0,\cds$,
\bea
\label{ulXn}
 \underline X_t^{t_0,\xi, n+1} = \xi+\int_{t_0}^t \wh b(\underline X_s^{t_0,\xi, n+1},\pa_x v(\cL_{\underline X^{t_0,\xi, n}}; s,\underline X_s^{t_0,\xi, n+1}),\cL_{\underline X_s^{t_0,\xi, n}})ds+B_t^{t_0}.
\eea
We then have the first main result of the paper.
\begin{thm}
\label{thm-minMFE}
Let Assumptions  \ref{assum-GH},   \ref{assum-b}, and \ref{assum-mon} hold. Then for any $(t_0, \mu)\in [0, T]\times \cP_2(\dbR^d)$ and $\xi\in \dbL^2(\cF_{t_0}; \mu)$, there exists a process $\ul X^{t_0, \xi}$ on $[t_0, T]$ such that\\
(i) $\ul X^{t_0, \xi, n}_t \preceq  \ul X^{t_0, \xi, n+1}_t$, $\forall n, t$, $\dbP$-a.s. with $\lim_{n\to\infty} \dbE[\sup_{t_0\le t\le T}|\ul X^{t_0, \xi, n}_t-\ul X^{t_0, \xi}_t|^2]=0$;\\
(ii) $\ul \nu^{t_0, \mu}:= \cL_{\underline X^{t_0,\xi}}$ is an MFE of the extended MFG at $(t_0, \mu)$;\\
(iii) for any MFE $\nu^*$ of  the extended MFG at $(t_0, \mu)$, we have $\ul \nu^{t_0,\mu} \preceq \nu^*$. That is, $\ul \nu^{t_0,\mu}$ is the minimal MFE.
\end{thm}
\proof For notational simplicity we omit $^{t_0, \xi}$ and $^{t_0,\mu}$.

First, by Assumption \ref{assum-b} and \reff{eq:paxvpaxxvbdd}, $\wh b(\underline X_s^{1},\pa_x v(\cL_{\underline X^{0}}; s,\underline X_s^{1}),\cL_{\underline X_s^{0}})\succeq - L^{\wh b}(C_1) {\bf 1}$. Then $\ul X^0_t \preceq \ul X^1_t$, $t_0\leq t\leq T$, $\dbP$-a.s. and thus $\cL_{\ul X^0} \preceq \cL_{\ul X^1}$. Applying Theorem \ref{thm:increase} repeatedly, we see that $\ul X^n$ is increasing in $n$, and thus we may define $\ul X := \lim_{n\to \infty} \ul X^n$. Moreover, following similar arguments one can easily see that $\ul X^n_t \preceq \ol X^0_t$, $t_0\leq t\leq T$, $\dbP$-a.s. for all $n$. Then it follows from the dominated convergence theorem that $\lim_{n\to\infty} \dbE[|\ul X^{n}_t-\ul X_t|^2]=0$, for any $t$. 

Next, by Assumption \ref{assum-b} and \reff{eq:paxvpaxxvbdd} we see that  $\h b(\cd, \pa_x v(\cd), \cd)$ is bounded by $L^{\h b}(C_1)$.  Then it follows from \cite[Lemma 4.1]{WZ} that the set $\cup_{n\ge 1}\{\cL_{\ul X^n_t}\}_{0\le t\le T}$ is precompact. Now send $n\to \infty$ in \reff{ulXn}, by the desired continuity of $\wh b$ in Assumption \ref{assum-b} and that of $\pa_x v$ in Lemma \ref{lem-vreg}, we have
\bea
\label{ulX}
 \underline X_t = \xi+\int_{t_0}^t \wh b(\underline X_s,\pa_x v(\cL_{\underline X}; s,\underline X_s),\cL_{\underline X_s})ds+B_t^{t_0}.
\eea
This implies that $\ul \nu:= \cL_{\underline X}$ is an MFE of the extended MFG at $(t_0, \mu)$. Moreover, compare this with \reff{ulXn}, one can easily see that $\lim_{n\to\infty} \dbE[\sup_{t_0\le t\le T}|\ul X^{n}_t-\ul X_t|^2]=0$.

Finally, for any MFE $\nu^*$ of  the extended MFG at $(t_0, \mu)$, consider the related SDE system:
\bea
\label{X*}
 X^*_t = \xi+\int_{t_0}^t \wh b(X^*_s,\pa_x v(\nu^*; s, X^*_s),\nu_s^*)ds+B_t^{t_0}.
\eea
Since $\nu^*$ is an MFE, we have $\nu^* = \cL_{X^*}$. Again since $\wh b(X^*_s,\pa_x v(\nu^*; s, X^*_s),\nu_s^*)\succeq - L^{\wh b}(C_1) {\bf 1}$, we have $\ul X^0_t \preceq X^*_t$, $t_0\leq t\leq T$, $\dbP$-a.s. Applying Theorem \ref{thm:increase} repeatedly, we see that $\ul X^n_t \preceq X^*_t$, $t_0\leq t\leq T$, $\dbP$-a.s. for all $n$. Then $\ul X_t\preceq X^*_t$, $t_0\leq t\leq T$, $\dbP$-a.s. and thus $\ul \nu \preceq \nu^*$.
\qed

We conclude this section with the following crucial flow property.
\begin{prop}
\label{prop-flow}
Let Assumptions  \ref{assum-GH},   \ref{assum-b}, and \ref{assum-mon} hold. Then, for any $(t_0, \mu)\in[0,T]\times\cP_2(\dbR^d)$,
\bea
\label{flow}
\ul\nu^{t_0, \mu}_t = \ul \nu^{t_1, \ul \nu^{t_0,\mu}_{t_1}}_t,\q \mbox{for all}~t_0\le  t_1\le t\le T.
\eea
\end{prop}
\proof Let $\xi\in \dbL^2(\cF_{t_0}; \mu)$. Then $\ul\nu^{t_0, \mu}_t  = \cL_{\ul X^{t_0,\xi}_t}$, $\forall t\ge t_0$, where $\ul X^{t_0,\xi}$ satisfies \reff{ulX}. Note that
\beaa
\underline X_t^{t_0,\xi} = \underline X_{t_1}^{t_0,\xi} +\int_{t_1}^t \wh b(\underline X_s^{t_0,\xi},\pa_x v(\cL_{\underline X^{t_0,\xi}}; s,\underline X_s^{t_0,\xi}),\cL_{\underline X_s^{t_0,\xi}})ds+B_t^{t_1},\q t\ge t_1.
\eeaa
We see that $\ul\nu^{t_0, \mu}$ is an MFE of the extended MFG at $(t_1, \cL_{\underline X_{t_1}^{t_0,\xi}}) = (t_1, \ul\nu^{t_0, \mu}_{t_1})$. Then by Theorem \ref{thm-minMFE} (iii) we have $ \ul \nu^{t_1, \ul \nu^{t_0,\mu}_{t_1}}_t\preceq \ul\nu^{t_0, \mu}_t$, for all $t\ge t_1$.

On the other hand, for the Picard iteration in \reff{ulX0} and \reff{ulXn},  by Theorem \ref{thm-minMFE} (i) we have $\ul X^{t_0,\xi, n}_{t_1} \preceq   \underline X_{t_1}^{t_0,\xi} =: \xi_1$, for all $n$. By \reff{ulX0} it is clear that $\ul X^{t_0,\xi,0}_{t} \preceq   \ul X^{t_1, \xi_1, 0}_t$ for all $t\ge t_1$. Note that
\beaa
 \underline X_t^{t_0,\xi, 1} = \ul X^{t_0,\xi, 0}_{t_1} +\int_{t_1}^t \wh b(\underline X_s^{t_0,\xi, 1},\pa_x v(\cL_{\underline X^{t_0,\xi, 0}}; s,\underline X_s^{t_0,\xi, 1}),\cL_{\underline X_s^{t_0,\xi, 0}})ds+B_t^{t_1}.
 \eeaa
Since  $\ul X^{t_0,\xi,1}_{t_1} \preceq  \xi_1$,  by Theorem \ref{thm:increase} we see that $\ul X^{t_0, \xi, 1}_t \preceq \ul X^{t_1,\xi_1,1}_t$, $t\ge t_1$, $\dbP$-a.s. Repeat the arguments, we obtain $\ul X^{t_0, \xi, n}_t \preceq \ul X^{t_1,\xi_1, n}_t$. Send $n\to\infty$, by  Theorem \ref{thm-minMFE} (i)  we have $\ul X^{t_0, \xi}_t \preceq \ul X^{t_1,\xi_1}_t$, $t\ge t_1$, $\dbP$-a.s. That is, $\ul\nu^{t_0, \mu}_t \preceq \ul \nu^{t_1, \ul \nu^{t_0,\mu}_{t_1}}_t$, for all $t\ge t_1$. Then we must have the equality.
\qed

\section{The corresponding value function}
\label{sec:value}
\setcounter{equation}{0}
In this section we investigate the dynamic value function corresponding to the minimal MFE:
\bea
\label{eq:minsol}
\underline V(t,x,\mu):=v( \ul \nu^{t,\mu}; t,x).
\eea
The following properties are immediate.
\begin{prop}\label{prop:paxVdecrease}
Let Assumptions  \ref{assum-GH},  \ref{assum-b}, and  \ref{assum-mon} hold. \\
(i) For any $(t,\mu)\in [0,T]\times\cP_2(\dbR^d)$, $\ul V(t, \cd, \mu) \in C^2(\dbR^d)$ with $|\pa_x \ul V|\le C_1$ and $|\pa_{xx} \ul V|\le C_2$ for the $C_1, C_2$ in \reff{eq:paxvpaxxvbdd};\\
(ii) for any $t\in [0, T]$, $\pa_x\underline V(t,\cd,\cdot)$ is increasing in $(x, \mu)$.
\end{prop}
\proof (i) is a direct consequence of Lemma \ref{lem-vreg} (i).

(ii) Assume $x_1 \preceq x_2$, $\mu_1\preceq \mu_2$ and let $\xi_i\in\dbL(\cF_{t_0},\mu_i)$, $i=1,2$, be such that $\xi_1\preceq \xi_2$. Then $\ul X_t^{t_0, \xi_1, 0} \preceq \ul X_t^{t_0, \xi_2, 0}$ for all $t_0\leq t\leq T$. Apply Theorem \ref{thm:increase} repeatedly, we have $\ul X_t^{t_0, \xi_1, n} \preceq \ul X_t^{t_0, \xi_2, n}$, $t_0\leq t\leq T$, $\dbP$-a.s. for all $n$. Then $\ul X_t^{t_0, \xi_1} \preceq \ul X_t^{t_0, \xi_2}$, $t_0\leq t\leq T$, $\dbP$-a.s. and hence $\ul \nu^{t_0,\mu_1} \preceq \ul \nu^{t_0,\mu_2}$.  Since $\pa_x\underline V(t,x,\mu)=\pa_xv( \ul \nu^{t,\mu}; t,x)$, then it follows from Proposition \ref{prop:umon} that $\pa_x\underline V(t,x_1,\mu_1) \preceq \pa_x\underline V(t,x_2,\mu_2)$.
\qed

However, as we will see in Section \ref{sect-eg} below, in general $\ul V$ is discontinuous in $(t, \mu)$. At below we show that $\pa_x \ul V$ is lower semi-continuous in $\mu$ in the following sense.

\begin{defn}
\label{defn-semicont}
(i) Let $\mu_n, \mu\in \cP_2(\mathbb R^d)$, $n\ge 1$. We say that $\mu_n\uparrow \mu$ (resp. $\mu_n \downarrow \mu$)  if $\mu_{n}\preceq \mbox{(resp. $\succeq$)} ~\mu_{n+1}$ for all $n$ and $\lim_{n\to\infty}W_2(\mu_n,\mu)=0$;\\
(ii) we say a function $U: \cP_2(\dbR^d)\to \dbR^d$ is lower  semi-continuous (resp. upper semi-continuous) if  $\dis\liminf_{n\to\infty} U(\mu_n) \succeq U(\mu)$ (resp. $\dis\limsup_{n\to\infty} U(\mu_n) \preceq U(\mu)$) whenever $\dis\lim_{n\to\infty} W_2(\mu_n, \mu)=0$. 
\end{defn}
Here $\liminf$ and $\limsup$ are taken component wise. We then have the semi-continuity of $\ul V$ in $(t,\mu)$.

\begin{prop}\label{prop:underlinecV}
Let Assumptions  \ref{assum-GH},  \ref{assum-b}, and  \ref{assum-mon} hold. Then \\
(i) for any $(t_k, \mu_k)\to (t, \mu)$, we have $\liminf_{n\to\infty} \pa_x \ul V(t_k, x, \mu_k) \succeq \pa_x \ul V(t,x,\mu)$, i.e. $\pa_x\ul V$ is lower semi-continuous in $(t,\mu)$.  Moreover, if $\mu_k \uparrow \mu$, then $\lim_{k\to\infty}\pa_x\underline V(t,x,\mu_k)  = \pa_x\underline V(t,x,\mu)$;\\
(ii) for any $x\in \dbR^d$ and $\nu\in C([0, T]; \cP_2(\dbR^d))$, the mapping $t\mapsto \pa_x \ul V(t, x, \nu_t)$ is lower semi-continuous, and in particular it is Borel measurable. 
\end{prop}
\proof  (i)  Fix $x$ and let $(t_k, \mu_k)\to (t, \mu)$, with $\xi_k\in \dbL^2(\cF_{t_k}; \mu_k)$, $\xi\in \dbL^2(\cF_t; \mu)$. Denote $\e_k := |t_k-t|+W_2(\mu_k, \mu)$ and $\hat t_k := t_k \vee t$. Then, by Proposition \ref{prop:umon} and \reff{paxvtreg} we have
\beaa
\left.\ba{c}
\dis \pa_x\underline V(t_k,x,\mu_k)=\pa_x v(\underline\nu^{t_k,\mu_k}; t_k, x)   \succeq \pa_x v(\cL_{\underline X^{t_k,\xi_k, n}}; t_k, x) \succeq \pa_x v(\cL_{\underline X^{t_k,\xi_k, n}_{[\hat t_k, T]}}; \hat t_k, x) - \rho(\e_k).
\ea\right.
\eeaa
 Recall \reff{ulX0} and \reff{ulXn}. It is clear that $\sup_{\hat t_k\le s\le T}W_2(\cL_{\ul X^{t_k, \xi_k, 0}_s},  \cL_{\ul X^{t, \xi, 0}_s}) \le \e_k+ \sqrt{\e_k}$. Similarly to the arguments in Theorem \ref{thm-minMFE}, we may utilize the locally uniform regularity in Assumption \ref{assum-b} with $R=C_1$ and with appropriate compact set $K$. Then, by  Lemma \ref{lem-vreg}  and stability of SDEs, one can easily show that there exists a modulus of continuity function $\rho_1$ such that $\sup_{\hat t_k\le s\le T}W_2(\cL_{\ul X^{t_k, \xi_k, 1}_s},  \cL_{\ul X^{t, \xi, 1}_s}) \le \rho_1(\e_k)$. Moreover, by  Lemma \ref{lem-vreg} and \reff{ulXn} again, we can show by induction on $n$ that  there exists a modulus of continuity function $\rho_n$ such that $\sup_{\hat t_k\le s\le T}W_2(\cL_{\ul X^{t_k, \xi_k, n}_s},  \cL_{\ul X^{t, \xi, n}_s}) \le \rho_n(\e_k)$. Then,  by  \reff{eq:abconpaxv}  and \reff{paxvtreg} we have, for each $n, k$,
\beaa
\dis \pa_x\underline V(t_k,x,\mu_k) \succeq \pa_x v(\cL_{\underline X^{t,\xi, n}_{[\hat t_k, T]}}; \hat t_k, x) -  \rho(\rho_n(\e_k))- \rho(\e_k) \succeq \pa_x v(\cL_{\underline X^{t,\xi, n}}; t, x) -  \rho(\rho_n(\e_k))- 2\rho(\e_k).
\eeaa
Send $k\to\infty$, we have $\liminf_{k\to\infty}\pa_x\underline V(t_k,x,\mu_k) \succeq \pa_x v(\cL_{\underline X^{t,\xi, n}}; t, x)$.
Now send $n\to\infty$, by \reff{eq:abconpaxv} again we have 
\beaa
\liminf_{k\to\infty}\pa_x\underline V(t_k,x,\mu_k) \succeq \pa_x v(\cL_{\underline X^{t,\xi}}; t, x) = \pa_x\underline V(t,x,\mu).
 \eeaa
 Moreover, if $\mu_k \uparrow \mu$, by Proposition \ref{prop:paxVdecrease} we have $\pa_x\underline V(t,x,\mu_k)  \preceq \pa_x\underline V(t,x,\mu)$, then the above inequality implies  $\lim_{k\to\infty}\pa_x\underline V(t,x,\mu_k)  = \pa_x\underline V(t,x,\mu)$.
 
 (ii) For $t_k \to t$, since $\nu_{t_k}\to \nu_t$, then $\liminf_{k\to\infty}\pa_x\underline V(t_k,x,\nu_{t_k}) \succeq  \pa_x\underline V(t,x,\nu_t)$. This proves the claimed lower semi-continuity, which implies further the Borel measurability. 
 \qed
 
\begin{defn}
\label{defn-cC2}
Let $\cC^2$ denote the set of functions $V: [0, T]\times \dbR^d \times \cP_2(\dbR^d)\to\dbR$ satisfying:\\
(i) $V(t,\cd, \mu)\in C^2(\dbR^d)$ for each $(t, \mu)$, and $\pa_x V, \pa_{xx} V$ are uniformly bounded;\\
(ii) for any $x\in \dbR^d$ and $\nu\in C([0, T]; \cP_2(\dbR^d))$, the mapping $t\mapsto \pa_xV(t,x, \nu_t)$ is Borel measurable.
\end{defn}
Then it is clear that $\ul V\in \cC^2$. The following lemma will be important in the next section.
\begin{lem}
\label{lem-SDEexistence}
Let Assumptions  \ref{assum-b} and \ref{assum-mon} (iii) hold and $V\in \cC^2$. Assume further that $\pa_x V$ is increasing in $\mu$ and lower or upper semi-continuous in $\mu$. Then, for any $(t_0, \mu)\in [0,T]\times\cP_2(\dbR^d)$ and $\xi\in \dbL^2(\cF_{t_0}; \mu)$, the following McKean-Vlasov SDE has a strong solution:
\bea
\label{X}
X_t^{t_0,\xi}=\xi+\int_{t_0}^t\wh b\big(X_s^{t_0,\xi},\pa_{x}V(s,X_s^{t_0,\xi}, \cL_{X_s^{t_0,\xi}}), \cL_{X_s^{t_0,\xi}}\big)ds+B_t^{t_0}.
\eea
Equivalently, the following Fokker-Planck equation has a weak solution $\nu(t,x)$:
\bea
\label{FP}
\pa_t \nu(t,x) - {1\over 2} \tr(\pa_{xx} \nu(t,x)) + div(\nu(t,x)\wh b(x, \pa_x V(t,x,\nu_t), \nu_t)) =0,\q \nu_{t_0} = \mu.
\eea
\end{lem}
\proof We shall only prove the case that $\pa_x V$ is lower semi-continuous in $\mu$. The upper semi-continuous case can be proved similarly, in the same spirit as we construct the maximal MFE in Subsection \ref{sect-max} below. 

Recall \reff{ulX0} and \reff{ulXn}. Denote $X^{t_0,\xi, 0} := \ul X^{t_0, \xi, 0}$, with possibly a larger $C_1$ which is an upper bound of $|\pa_x V|$, and for $n=0,1,\cds$,
\bea
\label{Xn}
X_t^{t_0,\xi, n+1} = \xi+\int_{t_0}^t \wh b(X_s^{t_0,\xi, n+1},\pa_x V(s, X_s^{t_0,\xi, n+1}, \cL_{X^{t_0,\xi, n}_s}),\cL_{X_s^{t_0,\xi, n}})ds+B_t^{t_0}.
 \eea
Since $\pa_x V$ is increasing in $\mu$ and by Assumption \ref{assum-mon} (iii), it is clear that $X^{t_0,\xi, n}$ is increasing in $n$, and $X_t^{t_0,\xi, n} \le \ol X_t^{t_0, \xi, 0}$ for all $t\in [t_0,T]$. Then there exists $X^{t_0, \xi}$ such that $\dis\lim_{n\to\infty} \sup_{t_0\le t\le T}\dbE[|X^{t_0,\xi, n}_t-X^{t_0,\xi}_t|^2]=0$. Note that, since $\pa_x V$ is increasing and lower semi-continuous in $\mu$, and $\cL_{X^{t_0,\xi, n}}\uparrow \cL_{X^{t_0,\xi}}$,  as in Proposition \ref{prop:underlinecV} (i) we have $\dis\lim_{n\to\infty}\pa_x V(t,x,  \cL_{X^{t_0,\xi, n}_t}) = \pa_x V(t,x,  \cL_{X^{t_0,\xi}_t})$. Then by sending $n\to\infty$ in \reff{Xn} we see that $X^{t_0,\xi}$ satisfies \reff{X}.
\qed

\section{Weak-viscosity solutions to the master equation}
\label{sect-minmax}
\setcounter{equation}{0}

\subsection{Viscosity solution to PDE system}
Differentiate \reff{HJB} formally in $x$, we  obtain the following system of PDEs: for $i=1,\cds, d$,
\bea
\label{eq:vecHJB} 
\pa_tu^i(t,x)+\frac{1}{2}\tr(\pa_{xx}u^i(t,x))+\pa_{x_i}H(x,u(t,x),\nu_t)+\pa_pH(x,u(t,x),\nu_t)\cd\pa_xu^i(t,x)=0.
\eea

\begin{defn}\label{def:viscosityHJBvec}
Fix $\nu\in C([0, T]\times \cP_2(\dbR^d))$ and consider $u: [0, T]\times \dbR^d\to \dbR^d$ such that both $u$ and $\pa_x u$ are bounded. We say that $u$ is a viscosity subsolution (resp. supersolution, solution) of the PDE system \reff{eq:vecHJB} if, for each $i$ and for given $u^{-i}:= (u^1, \cds, u^{i-1}, u^{i+1}, \cds, u^d)$, the function $u^i$ is a  viscosity subsolution  (resp. supersolution, solution) to the PDE  \reff{eq:vecHJB}  for fixed $i$ in the standard sense. 
  \end{defn} 
  
  \begin{lem}\label{lem:uviscosity}
Let Assumption \ref{assum-GH} hold true. Fix $\nu\in C([0,T];\cP_2(\dbR^d))$ and let $v(\nu; \cd,\cd)$ be the unique classical solution of the PDE \reff{HJB}. Then $u(t,x):=\pa_xv(\nu; t, x)$ is a viscosity solution to the PDE system \reff{eq:vecHJB}. 
\end{lem}
\proof Recall \reff{paxvrep} and \reff{tdY}. Note that $\td_x Y_s^{t,x,\nu} = u(s, X^{t,x}_s)$. Then, for fixed $i$, \reff{tdY} becomes:
\beaa
 \left.\ba{lll}
\dis \nabla_{x_i} Y_s^{t,x,\nu}=\pa_{x_i}G(X_T^{t,x},\nu_T)+\int_s^T[\pa_{x_i}H+\pa_pH\nabla_{x_i} Z_r^{t,x,\nu}](X_r^{t,x},  u^{-i}(r, X^{t,x}_r), \td_{x_i} Y_r^{t,x,\nu},\nu_r)dr\\
\dis \qq\qq\q -\int_s^T\nabla_{x_i} Z_r^{t,x,\nu}\cd dB_r.\\
\ea\right.
\eeaa
Then by the standard BSDE theory we see that $u^i(t,x) =  \nabla_{x_i} Y_t^{t,x,\nu}$ is a viscosity solution to the PDE \reff{eq:vecHJB} for each fixed $i$.
\qed

The next comparison principle is more or less standard, see e.g. \cite{HM} in slightly different contexts. We nevertheless sketch a proof for completeness. 
  \begin{lem}\label{lem:viscositycomparison}
Let Assumptions \ref{assum-GH} and  \ref{assum-mon} (i)-(ii) hold true, and fix $\nu\in C([0,T];\cP_2(\dbR^d))$. Let $u$ be as in Lemma \ref{lem:uviscosity}, and $\ul u$ and $\ol u$ be a viscosity subsolution and a viscosity supersolution, respectively, to the PDE system \reff{eq:vecHJB}  in the sense of Definition \ref{def:viscosityHJBvec}. 
If $\ul u(T,x)\preceq \pa_x G(x, \nu_T) \preceq  \ol u(T,x)$ for all $x\in \dbR^d$, then $\ul u\preceq u\preceq \ol u$ on $[0,T]\times\dbR^d$.
\end{lem}
\proof We shall prove only $\ul u \preceq  u$. The inequality $u\preceq \ol u$ can be proved similarly.

Fix $(t, x)$ and denote $X_s := x+ B^t_s$. For a possibly larger $C_1$ such that $ |\ul u|\le C_1$, recall \reff{tdY} and introduce the following linear BSDEs recursively: $\td_i Y^{0}:= C_1$, and for $n\ge 0$,
\bea
\label{tdYnu}
\left.\ba{c}
\dis \td_i Y^{n+1}_s = \pa_{x_i}G(X_T,\nu_T)-\int_s^T\nabla_i Z_r^{n+1}\cd dB_r\\
\dis +\int_s^T\Big[\pa_{x_i}H(X_r, \td^{-i} Y^n_r, \td_i Y^{n+1}_r, \nu_r)+\pa_pH(X_r, \td Y^n_r,\nu_r)\cd \nabla_{i} Z_r^{n+1}\Big]dr.
\ea\right.
\eea
That is, $\td Y^{n+1}_s = u_{n+1}(s, X_s)$, where $u^i_0 \equiv C_1$, and for $n\ge 0$ and for given $u_n$, the function  $u_{n+1}^i$ is the unique viscosity solution to the following PDE:
\bea
\label{un+1}
\left.\ba{c}
\dis \pa_tu^i_{n+1}(t,x)+\frac{1}{2}\tr(\pa_{xx}u^i_{n+1}(t,x))+\pa_{x_i}H(x,u_n^{-i}(t,x), u_{n+1}^i(t,x), \nu_t)\\
\dis +\pa_pH(x,u_n(t,x),\nu_t)\cd\pa_xu_{n+1}^i(t,x)=0,\q u^i_{n+1}(T,x) = \pa_{x_i} G(x, \nu_T).
\ea\right.
\eea
Recall \reff{tdY}. One can easily show that $ \lim_{n\to \infty} \sup_{t\le s\le T} \dbE[|\td Y^n_s - \td_x Y^{t,x,\nu}_s|^2] =0$, and thus $\lim_{n\to\infty} u_n = u$. Moreover, similar (actually easier) to the proof of Proposition \ref{prop:umon}, we can prove by induction on $n$ that $u_n$ is increasing in $x$ for all $n$. We claim that
\bea
\label{unu}
\ul u \preceq u_n,\q\mbox{for all}~n.
\eea
Then, by sending $n\to \infty$, we obtain $\ul u \preceq u$.

To see \reff{unu}, first, since $u_0^i\equiv C_1 \ge \ul u^i$, it holds true for $n=0$. Assume it holds true for $n$, and we shall verify it for $n+1$. By Assumption \ref{assum-mon} (ii) and $\pa_x u_{n+1}^i \succeq {\bf 0} $, we see that
\beaa
&&\pa_{x_i}H(x,u_n^{-i}(t,x), u_{n+1}^i(t,x), \nu_t)+\pa_pH(x,u_n(t,x),\nu_t)\cd\pa_xu_{n+1}^i(t,x)\\
&&\ge \pa_{x_i}H(x,\ul u^{-i}(t,x), u_{n+1}^i(t,x), \nu_t)+\pa_pH(x,\ul u(t,x),\nu_t)\cd\pa_xu_{n+1}^i(t,x).
\eeaa
Then $u^i_{n+1}$ is a viscosity supersolution of the following PDE:
\bea
\label{un+2}
\left.\ba{c}
\dis \pa_tu^i_{n+1}(t,x)+\frac{1}{2}\tr(\pa_{xx}u^i_{n+1}(t,x))+\pa_{x_i}H(x,\ul u^{-i}(t,x), u_{n+1}^i(t,x), \nu_t)\\
\dis +\pa_pH(x, \ul u(t,x),\nu_t)\cd\pa_xu_{n+1}^i(t,x)\le 0,\q u^i_{n+1}(T,x) = \pa_{x_i} G(x, \nu_T).
\ea\right.
\eea
Notice that $\ul u^i$ is a viscosity subsolution of the above PDE. Then by the standard comparison principle we obtain $\ul u^i \le u^i_{n+1}$. This proves \reff{unu} for $n+1$, and hence $\ul u \preceq u$.
\qed

\subsection{Weak-viscosity solutions to the master equation}

We now introduce a notion of weak-viscosity solution to the master equation \reff{master}, adapted from \cite{MZ2}. Recall Definition \ref{defn-cC2}. 
\begin{defn}
\label{def:viscosity}
We say that $V\in \cC^2$ is a weak-viscosity subsolution (resp. supersolution, solution) of the master equation \eqref{master} if, for any $(t_0, \mu)\in [0,T]\times\cP_2(\dbR^d)$, the Fokker-Planck equation \reff{FP} has a weak solution $\nu$ such that the function $u(t,x):=\pa_xV(t,x, \nu_t)$ is a viscosity subsolution (resp. supersolution, solution) to the PDE system \reff{eq:vecHJB}  on $[t_0, T]$ in the sense of Definition \ref{def:viscosityHJBvec} and satisfies $u(T,x)\preceq (\mbox{resp. $\succeq$, $=$}) \pa_x G(x,\nu_T)$.
\end{defn}

We first have the following simple result.
\begin{prop}
\label{prop-viscosity}
Let Assumptions  \ref{assum-GH}, \ref{assum-b}, and \ref{assum-mon} (i)-(ii) hold. Assume $V\in \cC^2$ is a weak-viscosity  solution of the master equation \eqref{master}. Then, for any $(t_0, \mu)\in [0,T]\times\cP_2(\dbR^d)$, the $\nu$ in Definition \ref{def:viscosity} is an MFE of the extended MFG at $(t_0, \mu)$.
\end{prop}
\proof First by Lemma \ref{lem-vreg} let $v(\nu;\cd,\cd)$ be the classical solution of the PDE \reff{HJB}. Then by Lemma \ref{lem:uviscosity} $\tilde u := \pa_x v(\nu;\cd,\cd)$ is a viscosity solution of the PDE system \reff{eq:vecHJB}  in the sense of Definition \ref{def:viscosityHJBvec} with $\tilde u(T, x) = \pa_x G(x, \nu_T)$. Now by Definition \ref{def:viscosity} and the comparison principle in Lemma \ref{lem:viscositycomparison}, we have $\pa_x v(\nu; t,x) = \pa_xV(t,x, \nu_t)$. This identifies \reff{FP} and \reff{Xb} with $\nu_t = \cL_{X^{t_0, \xi, \nu}_t}$, except that one is in PDE form while the other is in SDE form. Thus $\nu = \Phi(t_0, \mu, \nu)$, namely $\nu$ is an MFE at $(t_0, \mu)$.
\qed

\begin{rem}
\label{rem-viscosity}
Alternatively, we may call $V\in \cC^2$ a weak-viscosity  solution of the master equation \eqref{master} if it is both a weak-viscosity subsolution and a weak-viscosity supersolution of \eqref{master}, where the weak-viscosity subsolution and supersolution are defined in Definition \ref{def:viscosity}. Under this alternative definition, we may use one $\nu$ for the subsolution property and another different $\nu$ (and hence a different $u$)  for the  supersolution property. So this is weaker than Definition \ref{def:viscosity}, in particular, a weak-viscosity solution in this alternative sense does not necessarily provide an MFE as in Proposition \ref{prop-viscosity}. 
\end{rem}

\no Our second main result of the paper is the following.

\begin{thm}
\label{thm-ulV}
Let Assumptions  \ref{assum-GH}, \ref{assum-b}, and \ref{assum-mon} hold. \\
(i) $\underline V$ is a weak-viscosity solution to the master equation \eqref{master};\\
(ii) for any  weak-viscosity supersolution $V$ to the master equation \eqref{master}, we have
\bea
\label{ulVmin}
\pa_x\underline V ~\preceq~ \pa_x V. 
\eea
\end{thm}
\proof (i) Fix $(t_0, \mu)\in [0,T]\times\cP_2(\dbR^d)$. By Theorem \ref{thm-minMFE} and in particular \reff{ulX} we see that $\ul\nu^{t_0,\mu}$ is a weak solution to the Fokker-Planck equation \reff{FP} with $V=\ul V$. Moreover, by Proposition \ref{prop-flow} we have
\beaa
\ul u(t, x) :=\pa_x \ul V(t, x, \ul \nu_t^{t_0,\mu}) = \pa_x v(\ul \nu^{t, \ul \nu^{t_0,\mu}_t}; t,x) = \pa_x v(\ul \nu^{t_0,\mu}; t,x). 
\eeaa 
Then by Lemma \ref{lem:uviscosity} $\ul u$ is a viscosity  solution to the PDE system \reff{eq:vecHJB} with $\nu_t = \ul \nu^{t_0, \mu}_t$.  Moreover, $\ul u(T,x)=\pa_x G(x,\ul\nu_T^{t_0,\mu})$. Therefore, $ \ul V$ is a weak-viscosity solution to the master equation \eqref{master}.

 (ii) Let $V$ be an arbitrary weak-viscosity supersolution to the master equation \eqref{master}. For any $(t_0,\mu)\in[0,T]\times\cP_2(\dbR^d)$, let $\nu, u$ be as in Definition \ref{def:viscosity}. Then, for any $\xi\in \dbL^2(\mathcal{F}_{t_0};\mu)$, the McKean-Vlasov SDE \reff{X} has a strong solution  $X^{t_0,\xi}$ with $\nu=\cL_{X^{t_0,\xi}}$. Recall \reff{ulX0} and \reff{ulXn}. It is clear that $\ul X_t^{t_0,\xi, 0} \preceq X_t^{t_0,\xi}$ for all $t\in [t_0,T]$. Denote $\ul \nu^{t_0,\mu, 0} := \cL_{\ul X^{t_0,\xi, 0}} \preceq \nu$. Note that $\pa_x v(\ul \nu^{t_0,\mu, 0}; \cd,\cd)$ is a viscosity solution to the PDE system \reff{eq:vecHJB} with $\ul \nu^{t_0,\mu, 0}$ and by Proposition \ref{prop:umon}  $\pa_x v$ is increasing in $x$. Then by Assumption \ref{assum-mon} (ii) one can easily see that $\pa_x v(\ul \nu^{t_0,\mu, 0}; \cd,\cd)$ is a viscosity subsolution to the PDE system \reff{eq:vecHJB} with $\nu$. Moreover,  by Assumption \ref{assum-mon} (i), 
 \beaa
 \pa_x v(\ul \nu^{t_0,\mu, 0}; T, x) = \pa_x G(x, \ul \nu^{t_0,\mu, 0}_T) \preceq  \pa_x G(x, \nu_T) \preceq u(T,x), 
 \eeaa
 Since $u$ is a viscosity supersolution of this system, then by the comparison principle Lemma \ref{lem:viscositycomparison}, we have $\pa_x v(\ul \nu^{t_0,\mu, 0}; t, x)  \preceq  u(t,x) = \pa_x V(t,x,\nu_t)$ for all $(t,x)$. Denote
 \beaa
 \ul b(t,x) := \wh b(x, \pa_x v(\ul \nu^{t_0,\mu, 0}; t, x), \ul \nu^{t_0,\mu, 0}),\q b(t,x) :=  \wh b(x,  \pa_x V(t,x,\nu_t); \nu_t).
 \eeaa
 By Assumption \ref{assum-mon} (iii) one can easily see that $\ul b \preceq b$, and $\pa_{x_j} \ul b^i \ge 0$ for all $i\neq j$. Then, comparing \reff{ulXn} and \reff{X}, it follows from Lemma \ref{lem-SDEcomparison}  that $\ul X_t^{t_0,\xi, 1} \preceq X_t^{t_0,\xi}$, $t_0\leq t\leq T$, $\dbP$-a.s. Repeat the arguments we can show that $\ul X_t^{t_0,\xi, n} \preceq X_t^{t_0,\xi}$, $t_0\leq t\leq T$, $\dbP$-a.s. and $\pa_x v(\cL_{\ul X^{t_0,\xi, n}}; t, x)  \preceq  u(t,x)$ for all $n$. Send $n\to \infty$, by Theorem \ref{thm-minMFE} and Lemma \ref{lem-vreg} (ii) we see that $\ul X_t^{t_0,\xi} \preceq X_t^{t_0,\xi}$, $t_0\leq t\leq T$, $\dbP$-a.s. and $\pa_x v(\ul\nu^{t_0,\mu}; t, x)  \preceq  u(t,x)$. Therefore, $\pa_x \ul V(t_0,x,\mu) = \pa_x v(\ul\nu^{t_0,\mu}; t_0, x) \preceq u(t_0,x) = \pa_x V(t_0, x, \mu)$. Since $(t_0,x,\mu)$ is arbitrary, we conclude the proof.
 \qed

\section{Some extensions}
\label{sect-extension}
\setcounter{equation}{0}

\subsection{The maximal case}
\label{sect-max}

Similarly to Section \ref{sec:minMFE}, we can construct the maximal MFE as follows. Fix $(t_0,\mu)\in [0,T]\times \cP_2(\dbR^d)$ and $\xi\in \dbL^2(\cF_{t_0};\mu)$. Let $\ol X^{t_0,\xi, 0}$ be defined by \reff{ulX0}, and for $n\ge 0$,
\bea
\label{ulXnol}
 \ol X_t^{t_0,\xi, n+1} = \xi+\int_{t_0}^t \wh b(\ol X_s^{t_0,\xi, n+1},\pa_x v(\cL_{\ol X^{t_0,\xi, n}}; s,\ol X_s^{t_0,\xi, n+1}),\cL_{\ol X_s^{t_0,\xi, n}})ds+B_t^{t_0}.
\eea
Then, as in Theorem \ref{thm-minMFE} and Proposition \ref{prop-flow}, we have the following results.
\begin{thm}
\label{thm-maxMFE}
Let Assumptions  \ref{assum-GH},   \ref{assum-b}, and \ref{assum-mon} hold. Then for any $(t_0, \mu)\in [0, T]\times \cP_2(\dbR^d)$ and $\xi\in \dbL^2(\cF_{t_0}; \mu)$, there exists a process $\ol X^{t_0, \xi}$ on $[t_0, T]$ such that  \\
(i)  $\ol X^{t_0, \xi, n+1}_t  \preceq \ol X^{t_0, \xi, n}_t $, $\forall n, t$, $\dbP$-a.s. with $\lim_{n\to\infty} \dbE[\sup_{t_0\le t\le T}|\ol X^{t_0, \xi, n}_t-\ol X^{t_0, \xi}_t |^2]=0$;\\
(ii) $\ol \nu^{t_0, \mu}:= \cL_{\ol X^{t_0,\xi}}$ is an MFE of the extended MFG at $(t_0, \mu)$ and satisfies the flow property:
\bea
\label{flow2}
\ol\nu^{t_0, \mu}_t = \ol \nu^{t_1, \ol \nu^{t_0,\mu}_{t_1}}_t,\q \mbox{for all}~ t_0< t_1\le t\le T;
\eea
(iii) for any MFE $\nu^*$ of  the extended MFG at $(t_0, \mu)$, we have $\ol \nu^{t_0,\mu} \succeq \nu^*$. That is, $\ol \nu^{t_0,\mu}$ is the maximal MFE.
\end{thm}

We next define 
\bea
\label{olV}
\ol V(t,x,\mu):=v( \ol \nu^{t,\mu}; t,x).
\eea

\begin{thm}
\label{thm-olV}
Let Assumptions  \ref{assum-GH},  \ref{assum-b}, and  \ref{assum-mon} hold. \\
(i) $\ol V\in \cC^2$, $\pa_x \ol V$ is increasing in $(x,\mu)$ and upper semi-continuous in $(t,\mu)$. Moreover, if $\mu_k \downarrow \mu$, then $\lim_{k\to\infty}\pa_x\ol V(t,x,\mu_k)  = \pa_x\ol V(t,x,\mu)$; \\
(ii) $\ol V$ is a weak-viscosity solution to the master equation \eqref{master};\\
(iii) for any weak-viscosity subsolution $V$ to the master equation \eqref{master}, we have
\bea
\label{olVmax}
\pa_x V ~\preceq~ \pa_x\ol V . 
\eea
\end{thm}

\no The following result is an immediate consequence of  Theorems \ref{thm-ulV} and \ref{thm-olV}.
\begin{cor}
Let Assumptions  \ref{assum-GH},  \ref{assum-b}, and  \ref{assum-mon} hold. If $\underline V=\overline V$ on $[0,T]\times\dbR^d\times\cP_2(\dbR^d)$, then the master equation \eqref{master} admits a unique  weak-viscosity solution $V:=\underline V=\overline V$.
\end{cor}

\subsection{The decreasing case}
\label{sect-decrease}
In this subsection we replace Assumption \ref{assum-mon} with Assumption \ref{assum-mon2}.

\begin{thm}
\label{thm-decreasing}
Let Assumptions \ref{assum-GH},  \ref{assum-b}, and \ref{assum-mon2} hold true.  \\
(i)  $\pa_x v$ is decreasing in $(x, \nu)$, and $\Phi$ in increasing in $(\mu, \nu)$;\\
(ii) for any $(t_0, \mu)\in [0, T]\times \cP_2(\dbR^d)$, there exist MFEs $\ul \nu^{t_0, \mu}$ and $\ol \nu^{t_0, \mu}$ of the extended MFG at $(t_0, \mu)$ such that, for any other MFE $\nu^*$ of  the extended MFG at $(t_0, \mu)$, we have $\ul \nu^{t_0,\mu} \preceq \nu^* \preceq \ol \nu^{t_0,\mu}$;\\
(iii) the minimal MFE $\ul \nu^{t_0, \mu}$ and the maximal MFE $\ol \nu^{t_0,\mu}$ satisfy the flow property \reff{flow} and \reff{flow2}.
\end{thm}

Again we define the value functions: 
\bea
\label{ulolV}
\ul V(t,x,\mu):=v( \ul \nu^{t,\mu}; t,x),\q \ol V(t,x,\mu):=v( \ol \nu^{t,\mu}; t,x).
\eea

\begin{thm}
\label{thm-olV2}
Let Assumptions  \ref{assum-GH},  \ref{assum-b}, and  \ref{assum-mon2} hold. \\
(i) $\ul V, \ol V \in \cC^2$, $\pa_x\ul V$ is  decreasing in $(x,\mu)$ and upper semi-continuous in $(t,\mu)$, and $\pa_x\ol V$ is decreasing in $(x,\mu)$ and  lower semi-continuous in $(t,\mu)$;\\
(ii) $\ul V$, $\ol V$ are weak-viscosity solutions to the master equation \eqref{master};\\
(iii) for any weak-viscosity subsolution $V_1$ and  weak-viscosity supersolution $V_2$ to the master equation \eqref{master}, we have
\bea
\label{olVmax2}
\pa_x\ul V ~\succeq~ \pa_x V_1,\qq \pa_x \ol V ~\preceq~\pa_x V_2. 
\eea
(iv) If $\underline V=\overline V$ on $[0,T]\times\dbR^d\times\cP_2(\dbR^d)$, then the master equation \eqref{master} admits a unique weak-viscosity solution $V:=\underline V=\overline V$.
\end{thm}

\subsection{The  common noise case}
\label{sect-common}
In this subsection we study the extended mean field game with a common noise. We shall only consider the problem under  Assumption \ref{assum-mon}. The case under  Assumption \ref{assum-mon2} is similar.

Let $B^0$ be the common noise which is independent of $\dbF$, $\beta\geq 0$ be a constant and $\h\beta^2:=1+\beta^2$. For any $t_0\in [0,T]$, denote $B^{0,t_0}_t:=B_t^0-B_{t_0}^0$, $t\in [t_0,T]$ and $\dbF^{0,t_0}:=\{\mathcal{F}_t^{B^{0,t_0}}\}_{t_0\leq t\leq T}$. Let $C(\dbF^{0,t_0};\mathcal{P}_2(\dbR^d))$ denote the set of stochastic measure flow $\nu: [t_0, T]\times \O \to \cP_2(\dbR^d)$ which is $\dbF^{0, t_0}$-progressively measurable and continuous in $t$. Given any  $\nu\in C(\dbF^{0,t_0};\mathcal{P}_2(\dbR^d))$, consider the following backward stochastic PDE on $[t_0,T]$:
\bea
\label{BSPDE}
&&\dis d v(\nu; t, x) =  z(\nu; t,x)\cd dB_t^0 - \Big[\tr\big(\frac{\h\beta^2}{2} \pa_{xx} v(\nu;t,x) + \b\partial_x z^\top(\nu;t,x)\big) +H(x,\pa_x v(\nu;t,x),\nu_t)\Big]dt,\nonumber\\
&&\dis  v(\nu;T,x)= G(x, \nu_T),
\eea
where the solution pair $(v,z)$ is $\dbF^{0,t_0}$-progressively measurable. Given $\xi\in\dbL^2(\mathcal{F}_{t_0})$, we still use $X^{t_0,\xi,\nu}$ to denote the strong solution to the following SDE on $[t_0,T]$:
\bea
\label{Xbc}
X^{t_0, \xi, \nu}_t = \xi + \int_{t_0}^t \wh b\big(X_s^{t_0,\xi,\nu}, \pa_xv(\nu; s,X_s^{t_0,\xi,\nu}),\nu_s\big)ds+B_t^{t_0}+\beta B_t^{0,t_0}.
\eea
Introduce the Nash field $\Phi$ on $C(\dbF^{0,t_0};\mathcal{P}_2(\dbR^d))$: for any $(t_0, \mu)\in [0, T]\times \cP_2(\dbR^d)$ and $\xi\in \dbL^2(\cF_{t_0};\mu)$, 
\bea
\label{NashFieldc}
\dis \Phi(t_0, \mu, \nu) := \{\cL_{X^{t_0, \xi, \nu}_t|\mathcal{F}^{0, t_0}_t}\}_{t_0\le t\le T},\q\forall \nu\in C(\dbF^{0,t_0};\mathcal{P}_2(\dbR^d)).
\eea
Fix $(t_0, \mu)$, define MFE as a fixed point of $\Phi(t_0, \mu, \cd)$. Then the corresponding master equation becomes second order:
\bea
\label{masterc}
\left.\ba{c}
\dis \pa_t V + \frac{1}{2} \tr(\pa_{xx} V) +  H(x,\partial_x V,\mu) + \cM V =0, \q V(T,x,\mu) = G(x,\mu),\q \mbox{where}\\
\dis \cM V(t,x,\mu)  := \tr\Big( \int_{\dbR^d}\Big[\frac{\h\beta^2}{2} \pa_{\tilde x} \pa_\mu V(t,x, \mu, \tilde x)  + \pa_\mu V(t, x, \mu, \tilde x)\wh b^\top(\tilde x,\pa_x V(t, \tilde x, \mu),\mu)\\
\dis +\beta^2\pa_x\pa_\mu V(t,x,\mu,\tilde x)+\frac{\beta^2}{2}\int_{\dbR^d}\pa_{\mu\mu}V(t,x,\mu,\bar x,\tilde x)\mu(d\bar x)\Big]\mu(d\tilde x)\Big).
\ea\right.
\eea


\begin{thm}
\label{thm-common}
Let Assumptions \ref{assum-GH},  \ref{assum-b}, and \ref{assum-mon} hold true.  \\
(i)  $\pa_x v$ is increasing in $(x, \nu)$, and $\Phi$ in increasing in $(\mu, \nu)$;\\
(ii) for any $(t_0, \mu)\in [0, T]\times \cP_2(\dbR^d)$, there exist MFEs $\ul \nu^{t_0, \mu}$ and $\ol \nu^{t_0, \mu}$ of the extended MFG at $(t_0, \mu)$ such that $\ul \nu^{t_0,\mu} \preceq \nu^* \preceq \ol \nu^{t_0,\mu}$ for all other MFE $\nu^*$ of  the extended MFG at $(t_0, \mu)$;\\
(iii) the minimal MFE $\ul \nu^{t_0, \mu}$ and the maximal MFE $\ol \nu^{t_0,\mu}$ satisfy the flow property \reff{flow} and \reff{flow2}, respectively,  $\mathbb P$-a.s.
\end{thm}
Here, for any $\nu^i\in C(\dbF^{0,t_0};\mathcal{P}_2(\dbR^n))$, $i=1,2$, the partial order $\nu^1\preceq\nu^2$ is extended naturally:  $\nu^1_t\preceq\nu^2_t$ for all $t\in [t_0,T]$, a.s. The monotonicity of $\pa_x v$ in $(x, \nu)$ is also in obvious sense.

Define the value functions corresponding to the minimal and maximal MFEs respectively: 
\bea
\label{ulolVc}
\ul V(t,x,\mu):=v( \ul \nu^{t,\mu}; t,x),\q \ol V(t,x,\mu):=v( \ol \nu^{t,\mu}; t,x).
\eea
We note that $\ul V$ and $\ol V$ are $\cF^{0,t}_t$-measurable and hence are actually deterministic.
\begin{thm}
\label{thm-olV2c}
Let Assumptions  \ref{assum-GH},  \ref{assum-b}, and  \ref{assum-mon} hold. Then $\ul V, \ol V \in \cC^2$, $\pa_x\ul V$ is increasing in $(x,\mu)$ and lower semi-continuous in $(t,\mu)$, and $\pa_x\ol V$ is increasing in $(x,\mu)$ and upper semi-continuous in $(t,\mu)$.
\end{thm}

We may continue to study weak-viscosity solution of the master equation \reff{masterc} as in Section \ref{sect-minmax}. In this case the PDE \reff{eq:vecHJB} becomes a backward SPDE \reff{BSPDE}, which can be viewed as a path dependent PDE, see e.g. Zhang \cite[Chapter 11]{Zhang}.  However, in this case $\ul u(t,x, \o) := \pa_x V(t, x, \ul \nu^{t_0,\mu}_t(\o))$ is in general discontinuous in $(t, \o)$, thus the viscosity theory for path dependent PDEs in Ekren-Touzi-Zhang \cite{ETZ1,ETZ2} and Zhou \cite{Zhou} cannot be applied here. One possibility is to adapt the viscosity solution for backward SPDEs in Qiu \cite{Qiu}, which does not require the regularity in $\o$. On the other hand, we note that the value function \reff{Vextended} for the  MFG with a major player will have the same regularity issue, even when there is no common noise. So we shall leave the systematic investigation of this issue to a future research. 
%

\section{An example}
\label{sect-eg}
\setcounter{equation}{0}
In this section we solve an example completely. In particular, we shall show that $\ul V$ is in general discontinuous in $(t, \mu)$. Set $d=1$ and denote
\beaa
m(\mu) := \int_\dbR x\mu(dx).
\eeaa
Consider the example:
\bea
\label{eg-HGb}
 G(x,\mu):= xm(\mu), \q H(x,p,\mu):= \frac{p^2}{2},\q  \wh b(x, p, \mu):=\wh b(p):=\left\{\ba{lll}
\dis -2,& p<-2;\\
\dis 2p+{1\over 2} p^2,& -2\leq p<0;\ms\\
\dis 2p-{1\over 2} p^2,& 0\leq p<2;\\
\dis 2,& p\geq 2.
\ea\right.
\eea
One can easily verify that Assumptions  \ref{assum-GH}, \ref{assum-b}, and \ref{assum-mon} hold true.  Moreover,
\eqref{HJB} becomes
\beaa
\pa_t v(\nu; t,x )+\frac{1}{2}\tr(\pa_{xx}v(\nu; t,x))+\frac{|\pa_xv(\nu; t,x)|^2}{2}=0,\q v(\nu; T,x)=G(x,\nu_T) = x m(\nu_T).
\eeaa
It admits a unique solution:
\bea
\label{eg-v}
v(\nu; t,x)= x m(\nu_T) + {1\over 2}(T-t)|m(\nu_T)|^2.
\eea
Then $\pa_x v(\nu; t,x) = m(\nu_T)$ and thus \reff{Xb} becomes: 
\bea
\label{eg-Xb}
X^{t_0, \xi, \nu}_t = \xi + \wh b(m(\nu_T)) (t-t_0) +B_t^{t_0}.
\eea
Note that $\Phi$ depends on $\nu$ only through $m(\nu_T)$. Introduce the following operator:
\bea
\label{whPhi}
\wh \Phi(t_0, \mu, p) := \dbE\Big[\xi + \wh b(p) (T-t_0) +B_T^{t_0}\Big] = m(\mu) + \wh b(p) (T-t_0),\q p\in \dbR.
\eea
One can easily see that $\nu^*$ is an MFE at $(t_0, \mu)$ if and only if $p^* := m(\nu^*_T)$ is a fixed point of $\wh\Phi$: 
\bea
\label{p*}
p^* = \wh\Phi(t_0, \mu, p^*) = m(\mu) + \wh b(p^*) (T-t_0),\q\mbox{or equivalently,}\q \wh b(p^*) = {p^*-m(\mu)\over T-t_0}. 
\eea
Then, by \reff{eg-Xb} and \reff{eg-v}, the corresponding MFE and value are:
\bea
\label{eg-MFE}
X^{t_0, \xi, p^*}_t = \xi + \wh b(p^*) (t-t_0) +B_t^{t_0},\q v(p^*; t,x)= x p^* + {1\over 2}(T-t)|p^*|^2.
\eea

Note that one side of \reff{p*} is piecewise quadratic, and the other side is linear. By elementary calculation we  solve \reff{p*} in four cases. Denote
\bea
\label{eg-notation}
\left.\ba{c}
\dis \l := {1\over T-t_0},\q m_1 := 2-{2\over \l},\q m_2:= {\l\over 2} + {2\over \l}-2; \\
\dis \phi_-(\l, m) := \sqrt{(\l-2)^2 - 2\l m},\q \phi_+(\l, m) := \sqrt{(\l-2)^2 + 2\l m}
\ea\right.
\eea

{\it Case 1.} $\l\ge 2$. In this case $m_1>0$, and there is a unique fixed point $p^*$: 
\bea
\label{case1}
p^* = \left\{\ba{lll} m(\mu)-{2\over \l},\q \mbox{if}~ m(\mu) < -m_1;\ms\\
\l -2 - \phi_-(\l, m(\mu)),\q \mbox{if}~  -m_1\le m(\mu) <0;\ms\\ 
2-\l +  \phi_+(\l, m(\mu)), \q \mbox{if}~  0\le m(\mu) < m_1;\ms\\
m(\mu)+{2\over \l},\q \mbox{if}~ m(\mu) \ge m_1.
\ea\right.
\eea

{\it Case 2.} $4-2\sqrt{2}< \l < 2$.  In this case $0< m_2 < m_1$. We  solve the problem in three subcases.

{\it Case 2.1.} $|m(\mu)|> m_2$. In this case there is a unique fixed point: 
\bea
\label{case21}
p^* = \left\{\ba{lll} m(\mu)-{2\over \l},\q \mbox{if}~ m(\mu) < -m_1;\ms\\
\l -2 - \phi_-(\l, m(\mu)),\q \mbox{if}~  -m_1\le m(\mu) <-m_2;\ms\\ 
2-\l +  \phi_+(\l, m(\mu)), \q \mbox{if}~ m_2< m(\mu) < m_1;\ms\\
m(\mu)+{2\over \l},\q \mbox{if}~ m(\mu) \ge m_1.
\ea\right.
\eea

{\it Case 2.2.} $|m(\mu)| = m_2$. In this case  there are two fixed points:
\bea
\label{case22}
\left.\ba{c}
\dis p^* = \l -2 -  \phi_-(\l, m(\mu))\q \mbox{or}\q p^*= 2- \l,\q\mbox{if}\q m(\mu) =-m_2;\\
\dis p^* = \l -2\q\mbox{or}\q p^*=2-\l +  \phi_+(\l, m(\mu)), \q \mbox{if}\q m(\mu)=m_2.
\ea\right.
\eea

{\it Case 2.3.} $|m(\mu)| < m_2$. In this case there are three fixed points:
\bea
\label{case23}
\left.\ba{c}
\dis p^* = \l -2 -  \phi_-(\l, m(\mu))\q \mbox{or}\q p^*= 2-\l \pm  \phi_+(\l, m(\mu)),\q\mbox{if}\q -m_2<m(\mu) \le 0;\\
\dis p^* = \l -2 \pm  \phi_-(\l, m(\mu)) \q\mbox{or}\q p^*=2-\l +  \phi_+(\l, m(\mu)), \q \mbox{if}\q 0< m(\mu)<m_2.
\ea\right.
\eea

{\it Case 3.} $1< \l \le 4-2\sqrt{2}$.  In this case $0< m_1 \le m_2$. We  solve the problem in three subcases.

{\it Case 3.1.} $|m(\mu)|> m_2$. In this case there is a unique fixed point: 
\bea
\label{case31}
p^* = \left\{\ba{lll} m(\mu)-{2\over \l},\q \mbox{if}~ m(\mu) < -m_2;\ms\\
m(\mu)+{2\over \l},\q \mbox{if}~ m(\mu) > m_2.
\ea\right.
\eea

{\it Case 3.2.} $|m(\mu)| = m_2$. In this case  there are two fixed points:
\bea
\label{case32}
\left.\ba{c}
\dis p^* = m(\mu)-{2\over \l}\q \mbox{or}\q p^*= 2- \l,\q\mbox{if}\q m(\mu) =-m_2;\\
\dis p^* = \l -2\q\mbox{or}\q p^*=m(\mu)+{2\over \l}, \q \mbox{if}\q m(\mu)=m_2.
\ea\right.
\eea

{\it Case 3.3.} $|m(\mu)| < m_2$. In this case there are three fixed points:
\bea
\label{case33}
\left.\ba{c}
\dis p^* =m(\mu)-{2\over \l}\q \mbox{or}\q p^*= 2-\l \pm  \phi_+(\l, m(\mu)),\q\mbox{if}\q -m_2<m(\mu) \le -m_1;\\
\dis p^* = \l -2 -  \phi_-(\l, m(\mu))\q \mbox{or}\q p^*= 2-\l \pm  \phi_+(\l, m(\mu)),\q\mbox{if}\q -m_1<m(\mu) \le 0;\\
\dis p^* = \l -2 \pm \phi_-(\l, m(\mu)) \q\mbox{or}\q p^*=2-\l +  \phi_+(\l, m(\mu)), \q \mbox{if}\q 0< m(\mu)\le m_1;\\
\dis p^* = \l -2 \pm \phi_-(\l, m(\mu)) \q\mbox{or}\q p^*=m(\mu)+{2\over \l}, \q \mbox{if}\q m_1< m(\mu)< m_2;\\
\ea\right.
\eea

{\it Case 4.} $0< \l \le 1$.  In this case $0\le -m_1 < m_2$. We  solve the problem in three subcases.

{\it Case 4.1.} $|m(\mu)|> m_2$. In this case there is a unique fixed point: 
\bea
\label{case41}
p^* = \left\{\ba{lll} m(\mu)-{2\over \l},\q \mbox{if}~ m(\mu)  <-m_2;\ms\\ 
m(\mu)+{2\over \l},\q \mbox{if}~ m(\mu) > m_2.
\ea\right.
\eea

{\it Case 4.2.} $|m(\mu)| = m_2$. In this case  there are two fixed points:
\bea
\label{case42}
\left.\ba{c}
\dis p^* = m(\mu)-{2\over \l}\q \mbox{or}\q p^*= 2- \l,\q\mbox{if}\q m(\mu) =-m_2;\\
\dis p^* = \l -2\q\mbox{or}\q p^*=m(\mu)+{2\over \l}, \q \mbox{if}\q m(\mu)=m_2.
\ea\right.
\eea

{\it Case 4.3.} $|m(\mu)| < m_2$. In this case there are three fixed points:
\bea
\label{case43}
\left.\ba{c}
\dis p^* =m(\mu)-{2\over \l}\q \mbox{or}\q p^*= 2-\l \pm  \phi_+(\l, m(\mu)),\q\mbox{if}\q -m_2<m(\mu) \le m_1;\\
\dis p^* =m(\mu)\pm {2\over \l}, \q\mbox{or}\q p^*= 2-\l  -  \phi_+(\l, m(\mu)), \q \mbox{if}\q m_1\le m(\mu)<0;\\
\dis p^* =m(\mu)\pm {2\over \l}, \q\mbox{or}\q p^*= \l-2  +  \phi_-(\l, m(\mu)), \q \mbox{if}\q 0\le m(\mu)\le -m_1;\\
\dis p^* = \l -2 \pm  \phi_-(\l, m(\mu)) \q\mbox{or}\q p^*=m(\mu)+{2\over \l}, \q \mbox{if}\q -m_1< m(\mu)<m_2.
\ea\right.
\eea

Put all the cases together, we find that the minimal $p^*$, denoted as $\ul p^{t_0,\mu}$, is: 
\bea
\label{ulp*}
\ul p^{t_0,\mu} := \left\{\ba{lll}
\dis m(\mu) - {2\over \l},\q \mbox{if}\q \l>0, m(\mu) \le -m_1;\\
\dis \l -2 - \phi_-(\l, m(\mu)),\q \mbox{if}\q \l\ge 2, -m_1\le m(\mu) < 0,\\
\dis\qq\qq\qq\qq~\mbox{or}~ 0<\l<2, -m_1 < m(\mu) \le m_2;\\
\dis 2-\l +  \phi_+(\l, m(\mu)), \q \mbox{if}~\l\ge 2,  0\le m(\mu) < m_1,\\
\dis\qq\qq\qq\qq  ~\mbox{or}~ 4-2\sqrt{2} <\l<2, m_2 < m(\mu) < m_1;\\
\dis m(\mu)+{2\over \l},\q \mbox{if}~ \l> 4-2\sqrt{2}, m(\mu) \ge m_1,\\
\dis\qq\qq\qq\qq  ~\mbox{or}~ 0<\l\le 4-2\sqrt{2}, m(\mu) > m_2.
\ea\right.
\eea
By \reff{eg-MFE}, we then have the minimal MFE and the corresponding value function:
\bea
\label{eg-ulMFE}
\ul X^{t_0, \xi}_t = \xi + \wh b(\ul p^{t_0,\mu}) (t-t_0) +B_t^{t_0},\q \ul V(t_0,x, \mu)= x \ul p^{t_0,\mu} + {1\over 2}(T-t_0)|\ul p^{t_0,\mu}|^2.
\eea

Similarly, we find that the maximal $p^*$, denoted as $\ol p^{t_0,\mu}$, is: 
\bea
\label{olp*}
\ol p^{t_0,\mu} := \left\{\ba{lll}
\dis m(\mu) - {2\over \l},\q \mbox{if}\q \l>4-2\sqrt{2}, m(\mu) < -m_1,\\
\dis\qq\qq\qq\qq  ~\mbox{or}~ 0<\l\le 4-2\sqrt{2}, m(\mu) < -m_2;\\
\dis \l -2 - \phi_-(\l, m(\mu)),\q \mbox{if}\q \l\ge 2, -m_1\le m(\mu) < 0,\\
\dis\qq\qq\qq\qq~\mbox{or}~ 4-2\sqrt{2}<\l<2, -m_1 \le m(\mu) < -m_2;\\
\dis 2-\l +  \phi_+(\l, m(\mu)), \q \mbox{if}~\l\ge 2,  0\le m(\mu) < m_1,\\
\dis\qq\qq\qq\qq  ~\mbox{or}~ 0 <\l<2, -m_2 \le m(\mu) < m_1;\\
\dis m(\mu)+{2\over \l},\q \mbox{if}~ \l>0, m(\mu) \ge m_1;
\ea\right.
\eea
and the maximal MFE and the corresponding value function are:
\bea
\label{eg-olMFE}
\ol X^{t_0, \xi}_t = \xi + \wh b(\ol p^{t_0,\mu}) (t-t_0) +B_t^{t_0},\q \ol V(t_0,x, \mu)= x \ol p^{t_0,\mu} + {1\over 2}(T-t_0)|\ol p^{t_0,\mu}|^2.
\eea

 We note that, when $\l>2$, namely $T-t<{1\over 2}$,  $\ul p^{t, \mu}$ is smooth in $(t, \mu)$ and actually in this case $\ul V=\ol V$ is a classical solution of the master equation \reff{master}. This is consistent with the standard result that the master equation admits a unique classical solution over small time interval. 
 
 However, for $4-2\sqrt{2} < \l < 2$, namely ${1\over 2} < T-t<{1\over 4-2\sqrt{2}}$, we have
\bea
\label{discontinuous}
\left.\ba{c}
\dis \lim_{m(\mu) \uparrow m_2} \ul p^{t,\mu}  =  \l -2 - \phi_-(\l, m_2) = {1\over T-t} -2;\\
\dis \lim_{m(\mu) \downarrow m_2} \ul p^{t,\mu}  = 2-\l +  \phi_+(\l, m_2) = (1+\sqrt{2})(2-{1\over T-t});
\ea\right.
\eea
That is, $\pa_x \ul V(t,x,\mu) = \ul p^{t,\mu}$ is discontinuous in $\mu$ when ${1\over 2} < T-t<{1\over 4-2\sqrt{2}}$ and $m(\mu) = m_2$. 

Similarly, when $m(\mu) = {1\over 20}$, we see that $m_2 > m(\mu)$ if $T-t> {5\over 8}$ and $m_2  < m(\mu)$ if $T-t> {1\over 8}$. Then, by \reff{discontinuous} we have
\beaa
\lim_{t\uparrow (T-{5\over 8})} \ul p^{t,\mu}  = {1\over T-(T-{5\over 8})} -2=-{2\over5},\q \lim_{t\downarrow (T-{5\over 8})} \ul p^{t,\mu}=(1+\sqrt{2})(2-{1\over T-(T-{5\over 8})}) = {2(1+\sqrt{2})\over 5}.
\eeaa
That is, $\pa_x \ul V(t,x,\mu) = \ul p^{t,\mu}$ is discontinuous in $t$ at $t={1\over 8}$ and $m(\mu) ={1\over 20}$.

\end{document}